\documentclass[11pt]{article}
\usepackage{epsfig, psfrag, color, geometry, url}
\usepackage{graphicx}
\input epsf.sty
\epsfverbosetrue 

\usepackage{latexsym}
\usepackage{amsfonts}
\usepackage{amssymb}

\newcommand{\R}{\mathbb{R}}

\newcommand{\Z}{\mathbb{Z}}
\newcommand{\C}{\mathbb{C}}

\newcommand{\N}{\mathbb{N}}
\newcommand{\D}{\mathbb{D}}
\newcommand{\IS}{\mathbb{S}}

\newcommand{\rs}{\mbox{$\widehat{\C}$}}

\def\SSS{{\mathcal S}}

\def\DDD{{\mathcal D}}
\def\EEE{{\mathcal E}}
\def\FFF{{\mathcal F}}
\def\GGG{{\mathcal G}}

\def\MMM{{\mathcal M}}

\def\MMM{{\mathcal M}}
\def\WWW{{\mathcal W}}

\def\PPP{{\mathcal P}}
\def\RRR{{\mathcal R}}
\def\SSS{{\mathcal S}}
\def\UUU{{\mathcal U}}

\newtheorem{thm}{Theorem}[section]
\newtheorem{defn}[thm]{Definition}
\newtheorem{remark}[thm]{Remark}
\newtheorem{fact}[thm]{Fact}
\newtheorem{prop}[thm]{Proposition}

\newtheorem{lemma}[thm]{Lemma}
\newtheorem{cor}[thm]{Corollary}

\newcommand{\qed}{\nopagebreak \begin{flushright}
 \rule{2mm}{2.5mm} \end{flushright}}

\newcommand{\implies}{\Rightarrow}
\newcommand{\cl}{\overline}                      

\newcommand{\supp}{\mbox{\rm supp}}                    

\newcommand{\mod}{\mbox{\rm mod}}    
\newcommand{\intersect}{\cap}                    
\newcommand{\union}{\cup}                        
\newcommand{\diam}{\mbox{\rm diam}}					         
\newcommand{\dist}{\mbox{\rm dist}}
\newcommand{\lcm}{\mbox{\rm lcm}}      
\newcommand{\interior}{\mbox{int}}     

\newcommand{\mtwo}[4]                            
{\mbox{$\left(\begin{array}{cc}                  
#1 & #2 \\
#3 & #4 
\end{array}
\right)$}}

\newcommand{\dettwo}[4]                          
{\mbox{$\left|\begin{array}{cc}                  
#1 & #2 \\
#3 & #4 
\end{array}
\right|$}}

\newcommand{\pf}{\noindent {\bf Proof: }}

\newcommand{\be}{\begin{enumerate}}
\newcommand{\eb}{\end{enumerate}}
\newcommand{\bi}{\begin{itemize}}
\newcommand{\ib}{\end{itemize}}
\newcommand{\bl}{\begin{list}}
\newcommand{\lb}{\end{list}}

\newcommand{\gap}{\vspace{5pt}}                 

\def\dis{\displaystyle}
\newcommand{\al}{\alpha} 
\newcommand{\ep}{\varepsilon}
\newcommand{\g}{\gamma}
\newcommand{\G}{\Gamma}
\newcommand{\si}{\sigma}
\newcommand{\te}{\theta}
\newcommand{\de}{\delta}
\newcommand{\om}{\omega}

\newcommand{\cbar}{\rs} 
\newcommand{\roundness}{\mbox{\rm Round}}

\newcommand{\wtU}{{\widetilde{U}}}

\newcommand{\wtW}{\widetilde{W}}

\newcommand{\txi}{\tilde{\xi}}

\newcommand{\XX}{\mbox{$\mathfrak{X}$}}

\newcommand{\confdim}{\mbox{\rm confdim}}
\newcommand{\hdim}{\mbox{\rm H.dim}}
\begin{document}

\title{Minimal Ahlfors regular conformal dimension of coarse conformal dynamics on the sphere
}

\date{\today}

\author{Peter Ha\"{\i}ssinsky, Universit\'e de Provence \\ and \\ Kevin M.  
Pilgrim, Indiana University Bloomington}

\maketitle

\begin{abstract}
We prove that if the Ahlfors regular conformal dimension $Q$ of a topologically cxc map on the sphere $f: S^2 \to S^2$ 
is realized by some metric $d$ on $S^2$, then either $Q=2$ and $f$ is topologically conjugate to a semihyperbolic rational map 
with Julia set equal to the whole Riemann sphere, or $Q>2$ and $f$ is topologically conjugate to a map which lifts to  an affine expanding  
map of a torus whose differential has distinct real eigenvalues.   
This is an analog of a known result for Gromov hyperbolic groups 
with two-sphere boundary, and our methods apply to give a new proof.
\end{abstract}

\tableofcontents

\section{Introduction}

A topologically coarse expanding conformal (cxc) orientation-preserving branched covering map 
$f: S^2 \to S^2$ is an analog, in the setting of iterated maps, of a Gromov-hyperbolic group $G$ whose boundary 
is the two-sphere.  
Let  $\DDD$ denote the dynamical system on $S^2$ given either by the iterates of $f$ or by the elements of $G$.  
In both settings, there is a canonical quasisymmetry class $\GGG(\DDD)$, or conformal gauge, 
of Ahlfors-regular metrics on $S^2$ in which the elements of $\DDD$ have uniformly bounded metric distortion.  
With respect to a metric $d \in \GGG(\DDD)$, for maps, one says that $f$ is metrically cxc, while for groups, 
one says that the action of $G$ is uniformly quasi-M\"obius.     

When is $\DDD$ topologically conjugate to a genuine conformal dynamical system, i.e. to a rational map or 
to a Kleinian group? The topological assumptions imply that such a rational map is necessarily semihyperbolic and that such a Kleinian group acts cocompactly when extended to act on hyperbolic three-space. 
For groups, Cannon's Conjecture asserts that the answer  is always yes \cite{cannon:conj}. 
For maps, this is no longer the case:  
there are combinatorial obstructions to Euclidean conformality, discovered by Thurston \cite{DH1}.    

The Ahlfors regular conformal dimension $\confdim_{AR}(\DDD)$ of $\DDD$ is defined as the infimum of the set of 
Hausdorff dimensions $\hdim(S^2, d)$ of metrics $d \in \GGG(\DDD)$.  The general concept was introduced by Pansu, 
who computed a closely related  invariant for the natural metrics on the boundary at infinity 
of certain homogeneous manifolds of negative curvature \cite{pansu:confdim}. Since the topological dimension 
is always a lower bound, in our setting one has $\confdim_{AR}(\DDD) \geq 2$.   
If $\DDD$ is given by a group, then conjecturally, $\confdim_{AR}(\DDD)=2$.  
If $\DDD$ is given by a map, however, then the combinatorics of obstructions provide lower bounds on 
$\confdim_{AR}(\DDD)$ that may be strictly larger than $2$ \cite{kmp:ph:cxcIII}.  

A priori, the infimum in the definition may, or may not, be realized.  For groups with boundary $S^2$, 
Bonk and Kleiner \cite[Thm.\,1.1]{bonk:kleiner:conf_dim} have proven that if the infimum is achieved, 
then the conclusion of Cannon's Conjecture holds.   
For a general hyperbolic group, if the infimum is realized, then the group again has special properties; 
see \cite{kleiner:icm2006} for a survey of results.  
The main result of this work is an analogous statement in the setting of maps:  it gives a topological characterization of when the infimum is realized. 

\begin{thm}[Rational or Latt\`es]
\label{thm:rl} 
If $f: S^2 \to S^2$ is topologically cxc and a metric $d \in \GGG(f)$ realizes $\confdim_{AR}(f)$, 
then $f$ is topologically conjugate to either 
\be
\item a semihyperbolic rational map, in which case $\confdim_{AR}(f)=2$ , or 
\item an obstructed  Latt\`es example induced by an affine map on the torus whose differential has distinct positive real eigenvalues 
each larger than one, in which case $\confdim_{AR}(f)>2$.
\eb
\end{thm}

The Latt\`es examples lift to covering maps of tori.  
They are a ubiquitous family of exceptional cases to general statements in the dynamics of rational maps.  
Our theorem is additional evidence that they play a similar role for the dynamics of cxc maps; 
cf. \cite{martin:mayer:rigidity} and also the recent preprint \cite{qianyin:lattes}. 
Their conformal gauges are related to those arising from visual boundaries of certain three-dimensional solvable 
Lie groups.   

We will derive Theorem \ref{thm:rl} from the following general result, which applies to both maps and groups:

\begin{thm}
\label{thm:fol}
Suppose $\DDD$ is the dynamical system on $S^2$ determined by a topologically cxc map or a Gromov hyperbolic group.  
If $\confdim_{AR}(\DDD)$ is attained, then either 
\be
\item $\DDD$ is topologically conjugate to a semihyperbolic rational map or cocompact Kleinian group, or 
\item $\DDD$ preserves a foliation of $S^2$ having finitely many singularities.
\eb
\end{thm}

The two cases are not mutually exclusive; the overlap consists of so-called integral, or flexible, Latt\`es examples: 
the corresponding torus maps are conjugate to group endomorphisms of the form $x \mapsto m\cdot x$ where $m$ is an integer with $|m|>1$. 

  Theorem \ref{thm:fol} yields an alternative proof of the result of 
Bonk and Kleiner mentioned above:

\begin{thm}[M.\,Bonk \& B.\,Kleiner]\label{thm:hgcannon}
Let $G$ be a hyperbolic group with boundary homeomorphic to $S^2$.
Assume that its Ahlfors-regular conformal dimension is attained.
Then the action of $G$ on its boundary is  topologically
the action of a cocompact Kleinian group.
\end{thm}

The proof of Theorem \ref{thm:fol} relies on a new characterization of the gauge of the Euclidean two-sphere $\IS^2$ among those gauges on $S^2$ supporting such dynamical systems $\DDD$:  this gauge is the only such gauge containing a metric in which two rectifiable thick curves cross; 
see the Dichotomy in \S\,\ref{subscn:outlineproofs} below, \S\,\ref{subscn:analyticmoduli}, and Proposition \ref{prop:mainthmrat}.

In the remaining subsections of this introduction, we give precise definitions and outline the proofs.  

\subsection{Topologically cxc maps} 
Throughout this work, $f$ denotes a continuous, orientation-preserving, branched covering map from the sphere to itself with degree $\deg(f)\geq 2$. 
 The following class of dynamical systems was introduced in \cite{kmp:ph:cxci}.

\begin{defn}[Topologically cxc]
The map $f$ is called {\em topologically cxc} 
provided there exists a finite open covering $\UUU_0$ of $S^2$ by connected sets satisfying the following properties: 
\be  
\item[] {\bf [Exp]}  The mesh of the coverings $\UUU_n$ tends to zero as $n \to \infty$, 
where $\UUU_n$ denotes the set of connected components of $f^{-n}(U)$ as $U$ ranges over $\UUU_0$.  

\item[] {\bf [Irred]}  The map $f$ is {\em locally eventually onto}:  for any $x\in S^2$, and any neighborhood $W$ 
of $x$, there is  some $n$ with $f^n(W)=S^2$.

\item[] {\bf [Deg]} The set of degrees of maps of the form $f^k|\wtU: \wtU \to U$, where $U \in \UUU_n$, 
$ \wtU \in \UUU_{n+k}$, and $n$ and $k$ are arbitrary, has a finite maximum $p<\infty$.  
\eb
We denote by $\mathbf{U}=\union_{n \geq 0}\UUU_n$.  
\end{defn}
\gap

Note that the definition prohibits recurrent or periodic branch points.  
Also, one may choose $\UUU_0$ so that each atom of each cover $\UUU_n$ is
a Jordan domain.

A rational map is topologically cxc on the sphere if and only if it is semihyperbolic with Julia set the whole sphere, 
i.e. has neither nonrepelling cyles nor recurrent critical points \cite[Cor.\,4.4.2]{kmp:ph:cxci}.  
\gap

The {\em postcritical set} of $f$ is $P_f = \union_{n>0}f^n(B_f)$, where $B_f$ is finite 
the set of branch points of $f$.  If $P_f$ is finite, the {\em orbifold} associated to $f$ has weight function 
$\nu: S^2 \to \N$ given by $\nu(y) = \lcm\{\deg(f^n, x) : f^n(x)=y, n \geq 1\}$ where $\deg( , )$ is the local degree; cf. \cite{DH1}.

\subsection{Conformal gauges}
A homeomorphism $h: X \to Y$ between metric spaces is {\em quasisymmetric} provided 
there exists a homeomorphism $\eta: [0,\infty) \to [0,\infty)$ such that 
$d_X(x,a) \leq t\cdot d_X(x,b) \implies d_Y(f(x),f(a)) \leq \eta(t)\cdot d_Y(f(x),f(b))$ for all triples of points 
$x, a, b \in X$ and all $t \geq 0$.   

The conformal gauges associated to maps and to groups possess many metric and dynamical regularity properties.

A metric space $X$ is {\em Ahlfors regular}  of dimension
$Q$ provided there is a Radon measure $\mu$ such that
for any $x\in X$ and $r\in (0,\diam X]$,
$$\mu(B(x,r))\asymp r^Q\,.$$
If this is the case, this estimate also holds for the $Q$-dimensional Hausdorff measure.

\gap
The next concept is due to Mackay \cite{mackay:dim1}.  Suppose $L \geq 1$.  A metric space $X$ is {\em $L$-annularly linearly connected (ALC)} provided 
\bi
\item[](BT) for any $x$ and $y$ in $X$, there is a continuum $K$ containing both
points such that $\diam K\le L|x-y|$;
\item [](ALC2) for any $x\in X$, and $0<r\le 2r\le R\le\diam X$, any pair of points
in $B(x,R)\setminus B(x,r)$ can be joined by a continuum $K$ contained in
$B(x,LR)\setminus B(x,r/L)$.\ib
The initials (BT) stand for {\em bounded turning}.   
Such a metric space has no local cut points, and the continuum in condition (BT) can be taken to be an arc, 
indeed, a quasi-arc \cite{mackay:arcs}.  The portion of the next result dealing with maps summarizes \cite{kmp:ph:cxci}, Proposition 3.3.2 and Theorems 3.5.3 and 3.5.6; the portion dealing with groups may be found in \cite{paulin:determined} and \cite{mackay:dim1}.

\begin{thm}[Canonical Gauge]\label{thm:can_gauge}
Given a dynamical system $\DDD$ as above, 
we may endow $S^2$ with a distance $d_v$ and a measure $\mu$ with the following properties:
\be
\item the space $(S^2,d_v,\mu)$ is Ahlfors regular, annularly linearly locally connected and doubling;
\item the measure $\mu$ is quasi-invariant by $\DDD$, i.e. sets of measure zero are preserved;
\item 
\bi
\item[]$(\DDD=G)$ the action of $G$ is uniformly quasi-M\"obius: there exists an increasing homeomorphism
$\eta:\R_+\to\R_+$ such that, for all distincts quadruples $x_1,x_2,x_3,x_4$ and all $g\in G$,
$$[g(x_1):g(x_2):g(x_3):g(x_4)]\le \eta([x_1:x_2:x_3:x_4])$$
where $$[x_1:x_2:x_3:x_4] =\frac{d_v(x_1,x_2)d_v(x_3,x_4)}{d_v(x_1,x_3)d_v(x_2,x_4)}\,;$$
\item[] $(\DDD=\{f^n\})$
there exist constants $\te\in (0,1)$ and $r_0>0$ with the following properties:
$\diam_v U\asymp \te^n$ if $U\in\UUU_n$, $f^k(B_v(x,r\te^k))=B_v(f^k(x),r)$ for any $r<r_0$.\ib
\eb
Furthermore, if $d$ is another metric sharing these properties, then the identity map
between $(S^2,d_v)$ and $(S^2,d)$ is quasisymmetric.\end{thm}

It follows that the set $\GGG(\DDD)$ of all Ahlfors regular metric spaces $Y$ quasisymmetrically equivalent to $(S^2,d_v)$ 
is an invariant, called the {\em Ahlfors regular conformal gauge},  of the topological conjugacy class of $\DDD$; 
see \cite[Thm.\,3.5.3, Thm.\,3.5.6]{kmp:ph:cxci}.  
Therefore, the {\em Ahlfors regular conformal dimension}
\[ \confdim_{AR}(\DDD) := \inf_{Y \in \GGG(\DDD)} \hdim(Y)\]
is a numerical topological dynamical invariant as well.  
Moreover, this invariant almost characterizes  conformal dynamics among topological ones  on the sphere; 
see \cite[Thm.\,4.2.11]{kmp:ph:cxci} and \cite[Thm.\,1.1]{bonk:kleiner:qsparam}, \cite[Thm.\,1.1]{bonk:kleiner:rigidity}.

\begin{thm}[Characterization of conformal dynamical systems]
\label{thm:char_of_rm}
The dynamical system $\DDD$ is topologically conjugate to a
semihyperbolic rational map or cocompact Kleinian group  
if and only if $\confdim_{AR}(f) = 2$ and is realized.
\end{thm}
This is a consequence of Bonk and Kleiner's characterisation of the Riemann
sphere (Theorem \ref{thm:qsparam} in the Appendix below)
and of Sullivan's straightening of uniformly quasiconformal
groups and quasiregular maps.

\gap

\gap

\subsection{Real Latt\`es maps}
Let $\Gamma$ be the subgroup of isometries of the plane generated by $x \mapsto x+(1,0), x\mapsto x+(0,1)$, 
and $x \mapsto -x$.  Let $\Z^2 \subset \R^2$ be the standard integer lattice, and let $A$ be a 2-by-2 integer matrix. 

The quotient space $\R^2/\Gamma$ is homeomorphic to $S^2$.  The map 
$\tilde{f}_A: \R^2 \to \R^2$ given by $v \mapsto Av$ descends to a branched covering $f_A: S^2 \to S^2$.    
The collection of such maps $f_A$ generalizes the well-studied family of rational functions discovered by Latt\`es, 
so we call these {\em Latt\`es} maps.   
If the eigenvalues of $A$ are real, we say that $f_A$ is a {\em real Latt\`es map}.

The following facts are easily verified.  
A Latt\`es map $f_A$ is topologically cxc if and only if $A$ is expanding, i.e. all eigenvalues of $A$ lie 
outside the closed unit disk. The postcritical set $P_{f_A}$ of $f_A$ has at most four points, 
with equality if and only if the orbifold of $f_A$ has weight $2$ at each point of $P_{f_A}$.    
Conversely, any critically finite topologically cxc map $f: S^2 \to S^2$ whose orbifold has signature 
$(2,2,2,2)$ is topologically conjugate to some Latt\`es map $f_A$: the map $f$ lifts to an expanding map on 
the canonical double covering torus, and such toral maps are classified up to topological conjugacy by their action on homology.  Note that the linear map 
given by $A$ is determined only up to sign, cf. \cite[Prop.\,9.3]{DH1}. 

We will show 

\begin{thm}\label{thm:lattesclass}
Let $f_A$ be a Latt\`es map with associated  matrix $\pm A$.
\be
\item If $A$ is a multiple of the identity or has non-real eigenvalues, then 
$f_A$ is topologically conjugate to a   rational function, and $\confdim_{AR}(f)=2$ and is attained.
\item If $A$ is semisimple with real eigenvalues $1<|\lambda|<|\mu|$, then
$\confdim_{AR}(f_A)=1 +\log|\mu|/ \log|\lambda|$ and is attained.
\item If $A$ has a single repeated eigenvalue $|\lambda|>1$, then $\confdim_{AR}(f_A)=2$ and
is not attained.
\eb\end{thm}

The  above classification of Latt\`es maps is intimately
related to the classification of homogeneous three-manifolds of negative curvature.

Heintze \cite{heintze:homogeneous} has classified all homogeneous manifolds
of negative sectional curvature: they are solvable Lie groups obtained as an
extension of a nilpotent group by the group $\R$ associated to a derivation $\al$.
In \cite{pansu:confdim}, Pansu computes the conformal dimension of its boundary
at infinity when $\al$ is semi-simple, 
and gives a metric of minimal dimension. 
We complete his computations when $\al$ is not semi-simple in the three-dimensional case,
cf. \cite[\S\,6]{bonk:kleiner:conf_dim}:

\begin{thm}\label{thm:homoARdim}
Let $M$ be a homogeneous three-manifold of negative
sectional curvature given by a non-semi-simple derivation 
on $\R^2$. Then the Ahlfors-regular conformal dimension of its
boundary is $2$, but is not attained.\end{thm}

\subsection{Outline of proofs}\label{subscn:outlineproofs}

We begin by analyzing weak tangents $T$ of metric spheres $X=(S^2, d) \in \GGG(\DDD)$.  

We first prove a general result: any weak tangent $T$ of an arbitrary doubling and ALC metric surface is homeomorphic 
to the plane (Theorem \ref{thm:limit_surface}).   
In the particular case of $X \in \GGG(\DDD)$, the selfsimilarity of $X$ induced by the dynamics implies that $T$ is equipped with a locally 
quasisymmetric map $h: T \to X$; for groups, $h$ is a homeomorphism to the complement in $X$ of some point,  while for maps, it is a surjective branched covering map.  

Suppose the Ahlfors regular conformal dimension, $Q$, 
is attained by $X \in \GGG(\DDD)$. A theorem of Keith and Laakso \cite{keith:laakso} 
implies the existence of a weak tangent 
$T_1$ of $X$ which contains a family of curves of positive $Q$-modulus.  
By a theorem of Tyson, the associated map $h_1: T_1 \to X$ transports this family to such a family on $X$.    
Since $Q$-modulus behaves like an outer measure on the separable metric space of curves, 
there are ``density points'', 
called {\em thick curves}, on $X$.  The collection of thick curves on $X$ is invariant under the dynamics.

Two thick curves {\em cross} if any curve sufficiently close to  one intersects 
any curve sufficiently close to the other.  The key point in the proof of Theorem \ref{thm:fol} is the following 

\gap

\noindent{\bf Dichotomy.}  
\be
\item If there are thick curves which cross, then $Q=2$, $X$ is quasisymmetrically equivalent to the Euclidean round sphere, 
and the dynamics is conjugate to the action of a discrete cocompact group of M\"obius transformations
or to a rational map.
\item If there are no thick curves which cross, then there is a foliation $\FFF$ of $X$ by locally thick curves 
such that $\FFF$ has finitely many singularities and is invariant under the dynamics.
\eb

We prove (1) by means of a new estimate relating combinatorial  moduli in different dimensions (Proposition \ref{prop:comb_dim_estimate}
and Corollary \ref{cor:inter_surface_thick}) and a comparison theorem relating combinatorial and analytic moduli (Proposition \ref{prop:comb_mod}).  
To prove  (2), we show that  a tangent $T_2$ of $X$ at a suitable density point on a thick curve  
is foliated by curves whose images under the map $h_2: T_2 \to X$ are locally thick and which, by assumption, 
cannot cross; this yields the foliation $\FFF$ on $X$.  

For groups, case (2) cannot occur since the action of $G$ is minimal on $X$.
For maps, the invariance of $\FFF$ implies that $f$ is a Latt\`es map.
Theorem \ref{thm:lattesclass} concludes the proof of Theorem \ref{thm:rl}.

While the following is not used in the proof, it follows naturally from the arguments we give; it shows that the tangents arising in the proof can in fact be constructed directly.  

\begin{remark}
\label{rmk:orbifold-cover}
A posteriori we note that, in the case of maps, 
it is possible to construct from
a weak tangent at a fixed point the universal orbifold covering map $\pi:\tilde{X}\to X$ associated
to the Latt\`es map together with an expanding map  $\psi:\tilde{X}  \to \tilde{X}$ 
whose iterates are uniformly quasisymmetric and such that
$\pi\circ\psi=f\circ\pi$.
\end{remark}

\subsection{Outline of the paper} 
The core of the paper deals with the proof when $\DDD$ is a topologically cxc mapping, since this
is new. The case of hyperbolic groups is sketched in the appendix A3. 
In section 2, we relate Latt\`es maps to negatively curved homogeneous three-manifolds
and prove Theorems \ref{thm:lattesclass} and  \ref{thm:homoARdim}.
In section 3, we 
recall and complete results concerning weak tangents of Ahlfors-regular metrics in the gauge
of a topologically cxc map $f$ on $S^2$. We also prove (Theorem \ref{thm:limit_surface}) that weak tangent spaces of doubling annularly linearly connected 
surfaces are homeomorphic to the plane. Section 4 is devoted to the notion of moduli of curves and its coarse
version of so-called combinatorial moduli. There, we establish the new estimate
comparing combinatorial moduli in different dimensions. We conclude 
this section with the proof of statement (1) in the above dichotomy.
Section 5 is devoted to proving
that there families of curves of positive modulus, when the sphere is endowed with
a metric of minimal dimension. In section 6 we exhibit
a foliation by locally thick curves when no thick curve
cross: this proves statement (2) in the dichotomy and proves Theorem \ref{thm:fol} when $\DDD$ is topologically cxc.
In section 7, we summarize results and establish Theorem \ref{thm:rl}. 
In the appendix, we provide several applications of our estimate on combinatorial
moduli of different dimensions and we sketch the proof of Theorem \ref{thm:hgcannon} by first
establishing Theorem \ref{thm:fol}. We also state the results which still hold for general topologically cxc maps.

\subsection{Acknowledgements} We thank M. Bourdon, B. Kleiner, L. Kolev, and J. Los for useful conversations.  Both authors 
were partially supported by project  ANR Cannon (ANR-06-BLAN-0366). Part of this work was done while the first author
was visiting Indiana University. We thank its Mathematics Department for its hospitality.

\subsection{Notation and conventions}  

Throughout, $f$ denotes an orientation-preserving branched covering map of the two-sphere to itself of degree $d \geq 2$, unless otherwise specified. 
If $X$ is a metric space, the distance function $d$ on $X$ is understood, and $R>0$, 
we denote by $RX$ the metric space $(X, R\cdot d)$.  
Similarly, if $B=B(x,r)$ is a ball in $X$ and $c>0$, we set $cB = B(x, cr)$. When 
convenient we use the notation $|x-y|$ for $d(x,y)$.   For two positive functions $a, b$ we write $a \lesssim b$ 
if $a \leq C\cdot b$; $a \asymp b$ means $a \lesssim b$ and $b \lesssim a$.  
A sequence whose $n$th term is $a_n$ is denoted $(a_n)$. 
The round Euclidean two-sphere of constant curvature $+1$ is denoted $\IS^2$. 

\section{Classification of homogeneous manifolds and of Latt\`es maps}

In this section, we relate Latt\`es maps to negatively curved homogeneous three-manifolds and give a proof
of both Theorems \ref{thm:lattesclass} and  \ref{thm:homoARdim}.

\gap

Let $A:\R^2\to\R^2$ be an expanding linear map with integer coefficients.

We consider the Lie group $G$ obtained by the extension of the Abelian group
$\R^2$ by $A$, endowed with an invariant Riemannian metric. It follows from
Heintze \cite[Thm.\,1]{heintze:homogeneous} that $G$ is a negatively curved homogeneous
manifold. 

More precisely, we let $B$ be a matrix such that $\exp B= A$. We consider the three-dimensional Lie
algebra $\mathfrak{g}$ such that its derived algebra $\mathfrak{g}'=[\mathfrak{g},\mathfrak{g}]$ is Abelian
and two-dimensional; we fix a basis  $(v_1,v_2)$ of $\mathfrak{g}'$. We also assume that there exists a vector $v_0\in\mathfrak{g}$ such that,
for all $w\in\mathfrak{g}'$, $[v_0,w]=B(w)$, when written in the basis $(v_1,v_2)$. 
We  let $\langle\cdot,\cdot\rangle$ be the scalar product on $\mathfrak{g}$ which makes the basis $(v_0,v_1,v_2)$ orthonormal.
This defines a Riemannian metric on $G$.

Then the pointed boundary $\partial G\setminus\{\infty\}$ can be identified with $\mathfrak{g}'=\R^2$ 
endowed with a visual
metric. The Abelian subgroup $\R^2$ acts freely by isometries in this metric and 
the vector $v_0$ via the matrix $A$ induces a dilation $\tilde{f}_A$. 
The Lebesgue measure of $\R^2$ defines a so-called conformal measure on $\partial G\setminus\{\infty\}$.  

Taking a quotient of  $\partial G\setminus\{\infty\}$ by the group generated by $\Z^2$
and $-Id$ yields a metric two-sphere.  The dilation $\tilde{f}_A$ descends to a topologically cxc map---the
associated Latt\`es map $f_A$. By Theorem \ref{thm:can_gauge}, the metric belongs
to the canonical gauge of $f_A$,  so the Ahlfors-regular conformal dimension of $\partial G$ coincides with
the Ahlfors-regular conformal dimension of $f_A$. 

\gap 

\pf (Theorems \ref{thm:lattesclass} and  \ref{thm:homoARdim}) 
When $A$ is semisimple then the Ahlfors-regular conformal dimension of the boundary at infinity
of $G$ is attained \cite{pansu:confdim}: if there are two real eigenvalues $1<|\lambda| \le |\mu|$,
then $\confdim_{AR}\partial G= 1 +\log |\mu|/\log|\lambda|$. So if the eigenvalues are complex conjugate,
then $\confdim_{AR}\partial G=2$ and is attained. 
One may find an explicit metric on $\R^2$ of minimal dimension. In the real case,
let $(v_1,v_2)$ be a basis by eigenvectors associated to $\lambda$ and $\mu$ respectively.
Then set $$d(x_1 v_1 + x_2 v_2,y_1 v_1 + y_2 v_2)= |y_1-x_1|+|y_2-y_1|^{\al}$$
with $\al=\log|\lambda|/\log|\mu|$. One may easily check that for this metric, the map $\tilde{f}_A$ acts as a dilation,  
the corresponding Hausdorff dimension is $1+1/\al$, and the conformal dimension is minimized.
In the complex case, write $A= \lambda\cdot R$ where $R$ is a rotation. This defines canonically
a $\C$-linear map, so the  Ahlfors-regular conformal dimension is $2$.

We may now assume that $A$ is not semisimple. Then $A$ has a double eigenvalue
$e^\lambda$, $\lambda >0$. We may then find a basis in $\R^2$ such that
$$B=\mtwo{\lambda}{\de}{0}{\lambda}$$ with $\de$ arbitrarily small.
Using an explicit computation of the sectional curvatures given in the proof
of \cite[Thm.\,1]{heintze:homogeneous} yields that as $\delta \to 0$ the sectional curvature satisfies 
$K(x,y)=-\lambda^2 +O(\de)$ for any orthonormal vectors $x$, $y$.

It follows from \cite{bourdon:em} that one may choose the parameter $\ep$ of the visual distance to
be of order $\lambda + O(\de)$. Furthermore, the Bishop inequality \cite[Thm.\,3.101]{gallot:hulin:lafontaine}  implies that the volume entropy
$h(G)$ of $G$ is bounded from above by $2\lambda + O(\de)$. Finally, \cite[Lma\,5.2]{pansu:confdim} implies
that $\confdim_{AR}\partial G\le 2 +O(\de)$.
Letting $\de$ tend to $0$ establishes that $\confdim_{AR}\partial G=2$. 

It remains to prove that this dimension is not attained. Recall that $\tilde{f}_A$ defines a topologically cxc 
postcritically finite branched covering $f_A$ of the sphere with orbifold $(2,2,2,2)$. The visual
distances considered above all belong to the gauge of $f_A$ since $\tilde{f}_A$ acts as a dilation with
respect to this metric. If the conformal dimension were attained, then $f_A$ would be conjugate
to a rational map by Theorem \ref{thm:char_of_rm}. 
But, according to \cite[Prop.\,9.7]{DH1}, the matrix $A$ should then
be an integral multiple of the identity. Therefore, the conformal dimension is two but is not attained.
This ends the proof of both theorems.\qed

\section{Weak tangents}

We first recall the definitions of the Gromov-Hausdorff 
topology and of weak tangent
spaces. We then prove that weak tangents of ALC surfaces are homeomorphic to $\R^2$, and
we establish properties which will be used in the sequel.

\subsection{Gromov-Hausdorff convergence}

Let $Z$ be a proper metric space. Given two subsets $X,Y\subset Z$, define
their Hausdorff distance as
$$d_H(X,Y)=\max\left\{ \sup_{x\in X} d(x,Y),  \sup_{y\in Y} d(y,X)\right\}\,.$$

\gap

Note that $d_H(X,Y) = d_H(\cl{X}, \cl{Y})$.  
The set of nonempty compact subsets of $Z$ endowed with $d_H$ is complete 
and is sequentially compact if $Z$ is compact.

\gap

Given two compact metric spaces $X$ and $Y$,we  define their {\em Gromov-Hausdorff 
distance} by
$$d_{HG}(X,Y)= \inf_Z\{d_H(X,Y),\ X,Y\subset Z\}$$
where the infimum is taken over all compact metric spaces $Z$ which
contains isometric copies of both $X$ and $Y$.
This defines a distance function on the set of isometry classes of compact metric spaces.

\gap

In this topology, a sequence $(X_n)$ of compact metric spaces is sequentially relatively
compact if and only if the sequence is {\it uniformly compact} i.e., 
for any $\ep>0$, there exists an integer $N(\ep)$ such that
each set $X_n$ can be covered by $N(\ep)$ balls of radius $\ep$.

\gap

For a uniformly compact family of metric spaces $(X_\al)_{\al\in A}$, 
there exists a compact metric space $Z$ and isometric embeddings $X_\al\hookrightarrow Z$.
It follows that a sequence $(X_n)$ of compact metric spaces converges to $X$ with respect 
to the Gromov-Hausdorff distance if and only if there exists a fixed compact metric space $Z$ 
and isometric embeddings $j_n: X_n \hookrightarrow Z$ and $j: X \hookrightarrow Z$  
such that $j_n(X_n) \to j(X)$ in the Hausdorff topology on compact subsets of $Z$.  
Henceforth, we will regard the $X_n$ and the limit $X$ simply as subsets of $Z$.  
The completeness of the Gromov-Hausdorff metric implies the limit is in fact independent of $Z$.  

\gap

We now extend the notion of convergence to noncompact pointed spaces. 
Say a sequence of pointed proper metric spaces $(X_n,x_n)$ {\em converges in the Gromov-Hausdorff topology} 
to a pointed metric space $(X,x)$ if, for all radii $R>0$, 
the sequence of compact pointed closed balls $(\overline{B_{X_n}(x_n,R)},x_n)$ in $X_n$ 
converges with respect to the Gromov-Hausdorff distance to a compact set $B_R$ in $X$ which contains the open ball 
$(B_X(x,R),x)$.  

\gap

Finally, we extend the notion of convergence to subsets of pointed noncompact spaces.  
Suppose a sequence of pointed proper metric spaces $(X_n, x_n)$ converges.  
By definition, this implies that for all $R>0$, the sequence of compact closed balls $\overline{B_{X_n}(x_n, R)}$ 
converges.  This sequence must be uniformly compact, so there is a fixed compact metric space $Z_R$ 
into which each closed ball $\overline{B_{X_n}(x_n, r)}$ embeds.   
When $R_1 < R_2$ we may assume $Z_{R_1} \hookrightarrow Z_{R_2}$ isometrically and so we obtain an inductive system
of pointed proper metric spaces $(Z_R)_{R>0}$ containing embedded isometric copies of $(\overline{B(x_n,R)})_n$.

\begin{defn} Suppose $(X_n, x_n)$ is a convergent sequence of pointed proper metric spaces, 
and $Y_n \subset X_n$ are subsets.  The sequence $Y_n$ is said to {\em converge to $Y$} if 
\be
\item there exists $R_0>0$ such that $Y_n \intersect B_{X_n}(x_n, R_0) \neq \emptyset$ 
for all $n$ sufficiently large, and 
\item for all $R\geq R_0$, the compact sets  $\cl{Y_n} \intersect \cl{B_{X_n}(x_n, R)}$ 
converge in the Hausdorff topology on compact subsets of $Z_R$.  
\eb
\end{defn}

Convergence of maps is defined in the obvious way: suppose $(X_n, x_n) \to (X,x)$, $(Y_n, y_n) \to (Y,y)$, $f_n: (X_n, x_n) \to (Y_n,y_n)$, and 
$f: (X,x) \to (Y,y)$.  We say that $f_n \to f$ {\em in the Gromov-Hausdorff topology} if, after identifying spaces with their 
images under the above embeddings, we have $f_n(x_n) \to f(x)$ for every sequence $(x_n)$ 
with $x_n \in X_n$ for each $n \in \N$, such that $x_n \to x$.

\subsection{Weak tangents} 

The notion of weak tangent formalizes the processes of zooming in near a point and passing to a limit.  
This requires some tameness on the metric space.  A proper metric space $X$ with distance function $d$ is {\em $N$-doubling}, 
or {\em doubling} if $N$ bears no particular interest, if any
ball of finite radius can be covered by $N$ balls of half its radius.
In this case the family  $\{(X,x, Rd)\}_{x\in X,R>0}$ is relatively compact in the
Gromov-Hausdorff topology.  
A limit point of a sequence $(X, x_n, R_n d)$ of pointed rescaled spaces with $R_n\to\infty$ and $(x_n)$ lying in a compact
subset of $X$ is called a {\it weak tangent} of $X$. We speak of {\em weak tangents  at $x_0\in X$} if it is 
a limit of the sequence $(X,x_0, R_n d)$, i.e. $x_n=x_0$ is a constant sequence.

Note that any limit will be at most $2N$-doubling.

\gap

One may also consider limits of metric measure spaces $(X_n,x_n,\mu_n)$ where the $\mu_n$'s
are Radon measures. This means that $(X_n,x_n)$ converges as a sequence of metric spaces,
and in addition there exists constants $c_n$ such that $c_n\cdot\mu_n$ weakly converges.
In particular, if $(\mu_n)$ are all Ahlfors-regular measures of the same 
dimension
and with uniform constants, then, rescaling the measures if needed, 
one may extract a convergent subsequence to an Ahlfors regular measure $\mu$ in the limit.

\subsection{Limits of surfaces}

The main result of this subsection is the following:

\begin{thm}\label{thm:limit_surface}
Let $X$ be a doubling proper metric surface which is ALC.
Then any weak tangent space $T$ is ALC and is homeomorphic to the plane.
\end{thm}

The doubling property is used to ensure the existence of a tangent space, even if this
is not a necessary condition.

The proof has several steps. 
\be
\item Since the ALC condition is scale-independent, a straightforward argument shows $T$ is ALC.  
\item The definition of ALC implies immediately that the one-point, or Alexandroff, compactification $\hat{T}$ of $T$ is a locally connected continuum without local cut points.  
\item Next, we prove that $\hat{T}$ is embeddable in the sphere.  
To do this, we use the planarity of $X$ and a stability result regarding limits of graphs 
(Proposition \ref{prop:graph}) to show that $\hat{T}$ cannot contain certain non-planar graphs.  
We then appeal to a classical characterization theorem of Claytor \cite{claytor}.
\item  Therefore, the complement of $\hat{T}$ in the sphere is either empty, 
or consists of Jordan domains. The remainder of the proof consists of  ruling out this latter possibility.
We do this by analyzing the complement of simple closed curves in $T$.
\eb

Our graph stability result depends on three technical lemmas.

The first one says that tripods can be unzipped to a pair of arcs meeting at a point in such a way 
that the modification takes place only near two of the arms.  Think of $U$ as a small neighborhood of $\gamma'$. The arc $\alpha$ might enter and exit $U$ many times.

\begin{lemma}\label{lma:unzip}
Let $\g':[0,1]\to X$ be an arc in an $L$-ALC metric space
$X$ and let $U$ be   an arcwise connected open set which contains
$\g'([0,1))$. Let $\al':[0,1]\to X$ be an arc with the following
properties: $\al'(0)\notin U$, $\al'([0,1))\cap \g'=\emptyset$
and $\al'(1)\in (\g'\cap U)$. Then there are arcs
$\al: [0,1] \to X$ and $\g: [0,1] \to X$ such that 
\be
\item	\be
	\item $\gamma([0,1)) \subset U$,
	\item $\gamma(1)=\gamma'(1)$,
	\item there exists $0 < s_0 \leq 1$ such that $\gamma|_{[0, s_0]} = \gamma'|_{[0,s_0]}$;	
	\eb
\item there exists $0 < t_0 < 1$ such that 
\be
\item $\alpha|_{[0,t_0]} = \alpha'|_{[0, t_0]}$, 
\item $\alpha([t_0, 1)) \subset U$,
\item $\alpha([t_0, 1)) \intersect \gamma([0,1)) = \emptyset$,
\item $\alpha(1)=\alpha'(1)=\gamma'(1)=\gamma(1)$.  
\eb
\eb
\end{lemma}

\pf We will follow the 
same strategy as \cite[Lma\,3.1]{mackay:dim1}.
We will ``unzip'' the tripod
defined by $\g'$ and $\al'$ to have two curves which only meet at $x_0=\g'(1)$.

Let $x$ be the 
intersection point between $\g'$ and $\al'$. Let $r_x=\min\{|x-x_0|,d(x,
X\setminus U )\}/L$.
Pick a point $y\in \g'\cap ( B(x,r_x)\setminus B(x,r_x/2))$ beyond $x$, 
and a point 
$y'\in \al'\cap (B(x,r_x)\setminus B(x,r_x/2))$ before $x$.
By (ALC2), there is an arc $A$ in $U$ which joins $y'$ to $y$ in $B(x,Lr_x)\setminus  B(x,r_x/(2L))$.
If $A$ meets $\g'\cup \al'$ only at its extremities $y$ and $y'$, use $A$ to replace $\al'$ beyond 
$y'$. Otherwise, let $z$ be the first point in $\g'$ encountered 
beyond $x$ by $A$, and $z'$ the first
point before $z$ in $A$ which meets $\g'\cup\al'$: if $z'$ is on $\al'$, replace the piece of $\al'$ between
$z'$ and $x$ by $A$, which now ends at $z$; if $z'$ is on $\g'$, replace the part of $\al'$
between $z'$ and $z$ by $A$, and add to $\al'$ the piece of $\g'$ between $x$ and $z$.

This implies that we have now two curves $\g$ and $\al$ which meet at a definitely closer point from
$x_0$ than before, and they both continue with $\g'$ up to $x_0$; 
moreover,
they are contained in $U$. 
We continue inductively.  At each stage, modifications only alter the positions in $U$.  Since $r_x$ is locally bounded from below for $x\ne x_0$, 
both curves meet in the limit exactly at $x_0$.
\qed

The next lemma says that a collection of arcs whose endpoints are close can be extended to a multi-armed tripod.

\begin{lemma}\label{lma:join}
Let $U$ be an arcwise connected open set
of an annularly linearly connected metric space.
Let $n$ be a positive integer, and let us
consider $n$ pairwise disjoint arcs $\g_j'$, $j=1,\ldots,n$,
with $\g_j'(0)\notin U$, and $\g_j'(1)=x_j\in U$.
Let $x_0\in U$ be disjoint from these curves.
There are $n$ arcs  $\g_1,\ldots,\g_n$ such that $\g_j=\g_j'$
off of $U$ and which only meet at their other extremities
at $x_0$.
\end{lemma}

\pf We will procede by induction.
Let us assume that we have already constructed $\g_1,\ldots,\g_{k-1}$
such that the conclusion
of the lemma holds.

If $\g_k'$ intersects $\cup\g_j$, we let $\g$ be the one
which is first met. We may then apply Lemma \ref{lma:unzip}
to $\g$ and $\g_k'$ with the connected component of
$\{z\in U,\ d(z, \g) < d(z,(\cup\g_j)\setminus \g)\}$
which contains $\g\cap U$. This leads to $k$ arcs
which satisfy the conclusion of the lemma.

Otherwise, we consider an arc which joins $\g_k'(1)$ to $x_0$
in $U$. If it does not meet the requirements, we may apply
the argument above replacing $\g_k'$ with its extension
to the first point of intersection with the other curves.
\qed

Hausdorff convergence is quite weak.  If $X_n \to X$, arcs in $X$ need not be limits of arcs in $X_n$, even if the $X_n$ are required to be arcwise connected.
The lemma below says that if the geometry of the $X_n$ is controlled, then this is possible. 

\begin{lemma}\label{lma:arcs} 
Suppose $X$ is a compact subset of a proper metric space $Z$,  
$\g: [0,1] \to X$ is an arc, $\eta>0$, and $L\geq 1$.   
Then there exists $\varepsilon>0$ such that if $d_H(X,Y) \leq \varepsilon$ and $Y$ has $L$-bounded turning (BT), 
then, for all $y\in B(\g(0),\ep)\cap Y$ and $y'\in B(\g(1),\ep)\cap Y$
there exists an arc $c:[0,1]\to Y$ such that $c(0)=y$, $c(1)=y'$ and
$$\sup_{t\in[0,1]} |c(t)-\g(t)|\le \eta\,.$$
\end{lemma}

\pf Set $\ep=d_H(X,Y)$. We will construct an arc $c=c_\ep$ in $Y$ such that
$$\lim_{\ep\to 0} \sup_{t\in[0,1]} |c_\ep(t)-\g(t)|=0\,.$$

\gap

We first record some facts. Since $\g$ is an arc, there are increasing homeomorphisms
$\om_{\pm}:\R_+\to\R_+$ such that, for any $s,t\in [0,1]$,
$$\om_-(|s-t|)\le |\g(s)-\g(t)|\le \om_+(|s-t|)\,.$$

\gap

We note that, for any $0\le s < t\le 1$, one can find an arc $c:[s,t]\to Y$ such that
\begin{equation}\label{eq:1}
\sup_{u\in[s,t]} |c(u)-\g(u)|\le \ep + L(\om_+(|s-t|)+2\ep)+ \om_+(|s-t|)\,.\end{equation}

\gap

To prove this, pick $c(s),c(t)$ in $Y$ at distance at most $\ep$ from $\g(s),\g(t)$
respectively. By the (BT) property, there exists an arc $c:[s,t]\to Y$ such that
\begin{equation}\label{eq:0} \diam\, c([s,t])\le L|c(s)-c(t)|\le L(\om_+(|s-t|)+2\ep)\,.\end{equation}
Therefore, for any $u\in [s,t]$, the triangle inequality shows that
\begin{eqnarray*}
|c(u)-\g(u)| & \le & |c(u)-c(s)| + |c(s)- \g(s)| + |\g(s)-\g(u)|\\
& \le &  L(\om_+(|s-t|)+2\ep) + \ep + \om_+(|s-t|)\,.\end{eqnarray*}

\gap

Fix an integer $n\ge 1$ so that $\om_+(1/n)\le \ep$. For $j=0,\ldots n$, let $t_j=j/n$,
$x_j=\g(t_j)$, and consider $y_j\in Y$ such that $|x_j-y_j|\le \ep$ with $y_0=y$ and $y_n=y'$. By the procedure
above, build for $j=0,\ldots, (n-1)$, an arc $c_j:[t_j,t_{j+1}]\to Y$ joining $y_j$ to
$y_{j+1}$ such that, according to (\ref{eq:1}), for any $t\in [t_{j},t_{j+1}]$, 
\begin{equation}\label{eq:2}
|c_j(t)-\g(t)|  \le  (3 L+ 2)\ep \,.\end{equation}

\gap

Note that if $c_j\cap c_k\ne\emptyset$ for $j<k$, then, by (\ref{eq:0}),
$$|y_j- y_{k+1}|\le \diam\, c_j+\diam\,c_k \le 6L\ep\,,$$
so that $|x_j- x_{k+1}|\le (6L+2)\ep$. But since $ |x_j- x_{k+1}|\ge \om_-(t_{k+1} -t_j)$, it follows
that $$|t_{k+1} -t_j|\le \om_-^{-1} ((6L+2)\ep)\,.$$
Therefore, for any $t\in [t_j,t_{k+1}]$,
\begin{equation}\label{eq:3}
|\g(t)-x_j| \le (\om_+\circ \om_-^{-1})((6L+2)\ep)\,.\end{equation}

\gap

We extract an arc $c:[0,1]\to Y$ from the $c_j$'s by induction as follows. Let us first
define $\kappa:[0,1]\to Y$ as the concatenation of the arcs $c_j$'s. Set $s_0=0$ and $c(s_0)=\kappa(0)$.
If $s_{j-1}$ and $c|_{[0,s_{j-1}]}$ are constructed, let
$$u_j=\min\{t\in [s_{j-1},1],\ \exists\, s>t,\ \kappa(t)=\kappa(s)\}$$
if it exists or $u_j=1$ otherwise, 
and $$s_j=\max\{t\in [0,1],\ \kappa(t)=\kappa(u_j)\}\,.$$
If $u_j=1$, then we let $c|_{[s_{j-1},1]}=\kappa|_{[s_{j-1},1]}$.
Otherwise, set $c|_{[s_{j-1},u_j]}= \kappa|_{[s_{j-1},u_j]}$ and $c|_{[u_j,s_j]}=\kappa(u_j)$.
We may then continue. Since each $c_j$ is an arc, this procedure stops after at most
$n$ steps. We obtain a parametrized arc $c$ with connected fibers.

\gap 

We now estimate the sup-norm. If $t\in [s_j,u_{j+1}]$ for some $j$, then $c(t)$ coincides with $\kappa(t)$
so (\ref{eq:2}) implies
$$|c(t)-\g(t)|  \le  (3 L+ 2)\ep \,.$$

If $t\in [u_j,s_j]$, then $c(t)=\kappa (u_j)$ and 
there is some index $0\le k <n$ such that $t_k\le u_j <t_{k+1}$.
Applying (\ref{eq:0}) and (\ref{eq:3}), we obtain
\begin{eqnarray*}
|c(t)-\g(t)| & \le & |c_k(u_j)-y_k|+ |y_k- x_k|+ |x_k-\g(t)|\\
& \le & 3L\ep + \ep +(\om_+\circ \om_-^{-1})((3L+2)\ep)\,.
\end{eqnarray*}
\qed

\begin{prop}\label{prop:graph} Let $X$ be an $L$-annularly linearly  connected doubling proper metric space.
Let $T$ be a limit of $X_n=(X,x_n, R_n d)$ for $x_n\in X$ and $R_n\to \infty$.  
Suppose $f: \G \to \hat{T}$ is an embedding of $\G$.  
Then for all $n$ sufficiently large, there exists an embedding $f_n: \G \to X_n$.    
Moreover, if $f(\Gamma) \subset T$, then the $f_n$ may be chosen so that 
$\sup_{t \in \Gamma} |f_n(t)-f(t)| \to 0$ as $n \to \infty$.  
\end{prop}

\pf If $B$ is a ball in an ALC metric space, we let $B^0$ denote its connected component which
contains the center of the ball, so that $(1/L)B\subset B^0\subset B$.

Assume first that $f(\Gamma) \subset T$.  Then $f(\Gamma)$ is compact so we may replace the  spaces $X_n, T$ 
with closed balls, and we may assume the convergence of spaces is Hausdorff convergence 
in a compact metric space $Z$. 
It is convenient to supply $\G$ with a length distance which makes every edge isometric to $[0,1]$.
We may think of the restriction of $f$ to each edge $e$ as given by a map $f_e: [0,1] \to \Gamma$.  
Let $V$ be the vertex set of $\G$. 
Choose $r>0$ so that $|f(v)-f(w)| \geq 6Lr$ for each pair of distinct vertices $v, w$.

The continuity of $f$ implies the existence of $\de>0$ such that if $x\in\G$ and $d(x,V)\le \de$
then $d(f(x),f(V))\le r/2$. Furthermore, since $f$ is a homeomorphism onto its image, there is some
$\eta\in (0,r/2)$ such that, for any $(x,x')\in e\times e'$ with $e\ne e'$, if $|f(x)-f(x')|\le \eta$
then  $d(x,V)\le \de$. 
Since $X$ is ALC, the spaces $X_n$ have uniformly bounded turning. 
For each $n$ and each vertex $v\in V$, we pick points $x_v^n\in X_n$ such that $|f(v)-x_v^n|\le d_H(T,X_n)$.
According to Lemma \ref{lma:arcs}, there is some $n_0$ such that, for each $n\ge n_0$, for each edge $e$,
there is an arc $c^n_e: [0, 1] \to X_n$ approximating the edge $f_e$ with $c^n_e(v)= x_v^n$ and
$\sup|c^n_e(t)-f_e(t)| \leq \eta/2$.  
We let $f_n:\G\to X_n$ be the map 
obtained by setting $f_n|_e=c^n_e$, i.e. by concatenating each of the edge approximations. 

If $f_n$ is injective, then we are done. Otherwise, if $f_n(x)=f_n(x')$ for some $x\ne x'$, then $x$ and
$x'$ belong to different edges and $|f(x)-f(x')|\le |f(x)-f_n(x)|+|f_n(x')-f(x')|\le \eta$; therefore, $d(x,V)\le \de$,
$d(f(x),f(V))\le r/2$ and $d(f_n(x),x_v^n)\le \eta +r/2\le r$. Hence, $f_n$ is injective 
over $X_n\setminus (\cup_v B^0(x_v^n,Lr))$.
For each $v \in V$, 
let $e_1, \ldots, e_m$ be the edges incident to $v$, and assume that $f_{e_j}(0)=v$ for all $j$.
Apply Lemma \ref{lma:join} with $U=B^0(x^n_v, 2Lr)$ to truncated edges (so that they are disjoint)
to obtain an extension of the collection of arcs 
$\{c^n_{e_j}\}_{j=1}^m$ to arcs $[0,1] \to X_n$ meeting only at $x^n_v$.  
We have therefore produced an embedding 
$f_n: \Gamma \to X_n$. Note that if $f_n(x)\in B^0(x^n_v,2Lr)$, then $|f_n(x)-f(x)|\le 4 Lr +\eta/2\le 5Lr$.
Hence $f_n$ and $f$ are $5Lr$-close. 

Now suppose $f(\Gamma)$ meets the point at infinity.   
We may assume the point at infinity is the image of a vertex $v_\infty$.  
By bisecting the edges meeting the point at infinity we may assume there is a distinguished collection of valence 
two vertices $v_1, \ldots, v_m$ comprising the ends of a star-shaped graph with center at $v_\infty$.  
Suppose the edge $e_j$ joins the vertex $v_j$ to $v_\infty$, and the edge $g_j$ is the other edge incident to $v_j$, 
$j=1, \ldots, m$.  Let $\Gamma' $ be the subgraph obtained by deleting the edges $e_1, \ldots, e_m$, 
and let $\Gamma''$ be the subgraph obtained by deleting the edges $e_1, \ldots, e_m$ and the edges $g_1, \ldots, g_m$.  

By choosing $v_1, \ldots, v_m$  close enough to $v_\infty$, we may assume that $f(\Gamma'')$ is contained in a ball 
$B$ and that the images $f(v_1), \ldots, f(v_m)$ are contained in the ball $100LB$.  
By the first case, the restriction $f: \Gamma' \to T$ is approximated by $f_n: \Gamma' \to X_n$.
Let $B_n$ be a ball in $X_n$ which is close to $B$; let $U_n$ be the unbounded connected component of 
$X_n\setminus \cl{B_n}$.   
By the choice of $B$, the set $U_n$ is arcwise connected and contains the points 
$f_n(v_1), \ldots, f_n(v_n)$, which are endpoints of arcs $g^n_j$ approximating the arcs $g_j$.    
Choose arbitrarily a point $x_n^\infty \in U_n$ which does not meet the image of any arc $g_j^n$.  
By Lemma \ref{lma:join}, there is an extension of the arcs $g_j^n$ to a collection of arcs 
which meet only at $x_n^\infty$.  
We have produced an embedding of $\Gamma$ into $X_n$ and the proof is complete.\qed

\begin{remark}
The proof shows that if $Z_n$ are compact and uniformly ALC and if $Z_n \to Z$ in the Gromov-Hausdorff topology, 
then any finite graph embedding $f: \Gamma \to Z$ is a limit of embeddings $f_n: \Gamma \to Z_n$.
\end{remark}

We now prove the theorem.  We first record some facts and notation.
Let $X$ be an unbounded ALC proper metric space. Given any compact set
$K$ and a point $x\in K$, there is a single
connected component of $X\setminus K$ which contains points from $X\setminus B(x, 2L\diam K)$.
Therefore the {\it filling-in} of $K$ consisting of $K$ with the other components of $X\setminus K$
is a compact subset of $B(x,2L\diam K)$.

\gap

\pf (Thm\,\ref{thm:limit_surface}). Let $(R_n)$ tend to $+\infty$ and
$(T,t,d_T)$ be a limit of $X_n =(X,x,R_n d)$.
We may assume that the spaces are homeomorphic to planes for the statement is local.
Since the Aleksandroff compactification $\hat{T}$ of $T$ is locally connected with no
local 
cut points,  Claytor's imbedding theorem implies
that it is embeddable in the sphere if
$\hat{T}$ has no embbeded complete graphs on five vertices, nor a complete bipartite
graph on three vertices \cite{claytor}. This is the content of Proposition \ref{prop:graph}: if
there were such a graph, then $X_n$ would also have a copy in a Jordan neighborhood of $x$,
which is impossible.

This implies that we may think of $\hat{T}$ as a locally connected continuum of the sphere with
no local cut points. If this is the sphere then the proof is complete. Otherwise,
each connected component of the complement is a Jordan domain, see \cite[Thm.\,VI.4.4]{whyburn:analytic_topology}.

\gap

{\noindent\sc Claim.---} {\it For any simply closed curve $\g$ in $T$, there are two continua of $\hat{T}$, 
$\Omega$
and $D$, with $\Omega \intersect T$ unbounded in $T$ and $D$ bounded,   
such that $\Omega\cap D=\g$, $\Omega\cup D=\hat{T}$ and, for any $z\in \g$
and any $r>0$, both intersections $(\Omega\setminus\g)\cap B(z,r)$ and $(D\setminus\g)\cap B(z,r)$
are nonempty.}
\gap

The claim implies that $S^2\setminus \hat{T}$ has no bounded components.
Since the point at infinity is not a local cut point of $\hat{T}$, there can
be at most one unbounded component, so that $\hat{T}$ is a closed disk. 
Applying the Claim to
a bounded Jordan curve $\gamma$ in $T$ which follows a piece of its boundary, the claim also leads
us to a contradiction. So, up to the claim, the theorem is proven.

\gap

{\noindent\sc Proof of the claim. ---}  
Let $\g$ be a simple closed curve in $T$.   Choosing arbitrarily a pair of vertices, we may regard it
as a graph.  Proposition \ref{prop:graph} implies that
$\g$ is a limit of parametrized Jordan curves $\g_n$ in $X_n$ where the parametrizations converge uniformly. 
Each $\g_n$ for $n$ large enough
separates $X_n$ into two components $D_n$ and $\Omega_n$, the latter being unbounded. 
We may assume that $\cl{D}_n \to D$ and $\cl{\Omega}_n \to \Omega$ where 
$D, \Omega$ are closed connected subsets of $T$. Adding the point at infinity to $\Omega$ 
provides us with a covering of $\hat{T}$ by two continua.

In this paragraph, we prove that $\gamma = \Omega \intersect D$.  Suppose first that
$z\in\Omega\cap D$. Then by definition, there exist $w_n \in D_n$, $w_n' \in \Omega_n$ with $w_n \to z$ 
and $w_n' \to z$ as $n \to \infty$.   
The (BT) property implies that there is a curve $c_n$ joining 
$w_n$ to $w_n'$ such that $\diam\, c_n\le L \cdot d_n(w_n,w_n')$.
The Jordan curve theorem implies that $c_n$ intersects $\g_n$ in a point $z_n$.  
As $n$ tends to infinity, $z_n \to z$, so we obtain that $z\in\g$.  
Conversely, suppose now $z \in \gamma$.  
Then $z = \lim z_n$ with $z_n \in \gamma_n$.  
Applying the Jordan curve theorem  again gives $w_n \in D_n$ and $w_n' \in \Omega_n$ arbitrarily close to $z_n$; 
we may therefore assume $w_n \to z$ and $w_n' \to z$ and so $z \in \Omega \intersect D$. 

For the second part of the claim, we argue by contradiction.  
Suppose one of the intersections in the statement is empty;
we will only treat the case of $D$. So, assume that
there is some $z\in\g$ and $r>0$ such that $D\cap B(z,r)\subset \g$.     
Fix $0 < \ep \ll r$;
we will prove that this implies the existence of points $z_n,z_n'\in\g_n$ such that
the connected components of $\g_n\setminus\{z_n,z_n'\}$ have diameter comparable to
that of $\g_n$ but that $d_n(z_n,z_n')\le 2\ep$.  Letting $\ep$ go to $0$ and 
going to the limit, this will
contradict that $\g$ is a simple closed curve. 

Fix a point $y$ on $\gamma$ at distance $\frac{1}{2}\diam\, \g$ from $z$. Suppose $y_n, z_n \in \g_n$ and $y_n \to y$ and $z_n \to z$.   By assumption, for $n$ large enough,
$D_n\cap B(z_n,r/2)$ is contained in the $\ep$-neighborhood of $\g_n$. 

Let $E_n$ be the closure of the filling-in of the component of $B(z_n,r/6L^2)$
containing $z_n$ and $F_n$ be the closure of the unbounded component of $X_n\setminus B(z_n,r/2L)$
so that $d(E_n,F_n)\ge (r/3L)$ and $A_n \subset B(z_n, r/2)$ 
hold.   We may assume that $r$ is small enough so that $F_n$ contains $y_n$.
The set $A_n=X_n\setminus (E_n\cup F_n)$ is an annulus, and the curves $\gamma_n$ must connect the ends (and may also connect an end with itself); see Figure 1.  
\begin{figure}
\begin{center}
\includegraphics[width=2.5in]{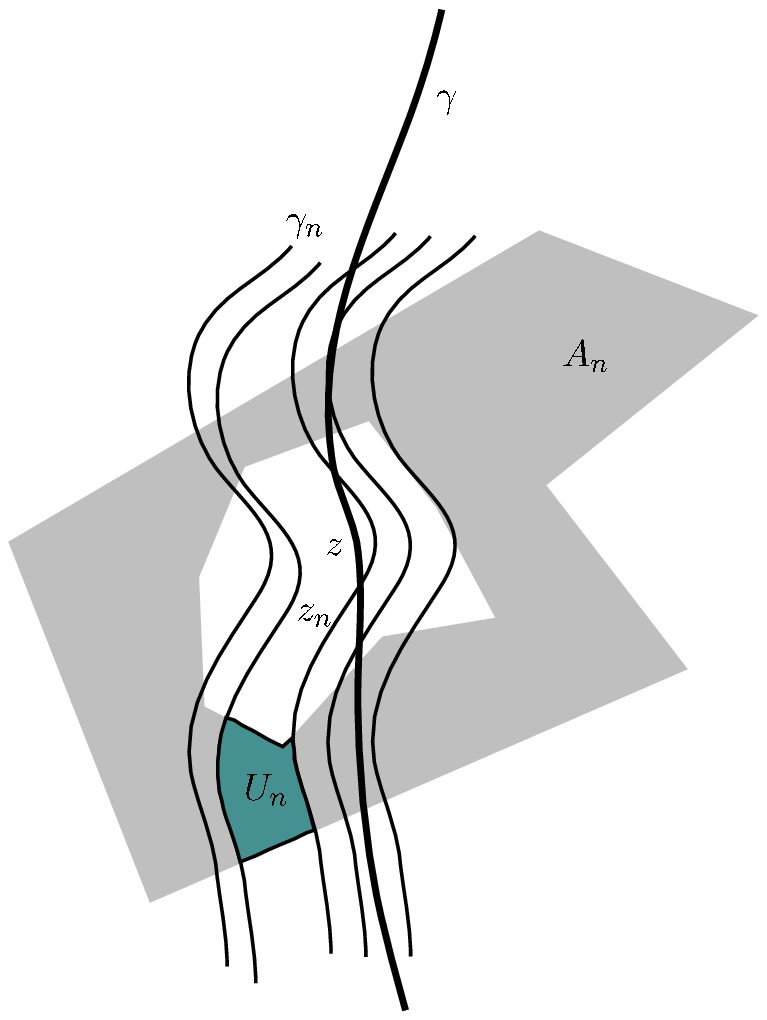}\\
{\sf Figure 1}
\end{center}

\end{figure}
Denote by $c_1^n$ and $c_2^n$
the {closures of the two components
of $\g_n\setminus\{z_n,y_n\}$; they are given by parametrizations (which we denote also by $c_j^n$) defined on compact intervals.  By the uniform continuity of the $c_j^n$, there is some $\de>0$ such that
if $|s-t|\le \de$ then $|c_j^n(s)-c_j^n(t)|\le d(E_n,F_n)/2$.  Since the $c_j^n$ are defined on compact intervals, it follows that there are 
only finitely many components of the intersection $c_j^n\cap A_n$ joining the end $E_n$ to the end $F_n$. 
Hence there is a component $U_n$ of $A_n\cap D_n$ which has in its
boundary pieces of both $c_1^n$ and of $c_2^n$.   
By assumption, $D_n \intersect B(z_n, r/2)$ is contained in 
the $\ep$-neighborhood of $\g_n$, and $A_n \subset B(z_n, r/2)$, 
so the component $U_n$ of the intersection is contained in the $\ep$-neighborhood of $\g_n$.

It follows that for all $u \in U_n$, we have that $d(u,\partial U_n) \le d(u, \gamma_n \intersect \partial U_n) \le  \ep$ where $d(u,\partial U_n)$ is the minimum distance from $u$ to a point in $\partial U$.   For each $j=1,2$, the set $U_n^j$ of points in $U_n$ which
are at distance strictly less than $\ep$ from $c_j^n$ is non-empty and open, 
hence the connectivity
of $U_n$ implies that at least one point in $U_n$ is at distance at most $\ep$ from 
both curves $c_1^n$ and $c_2^n$: we may find $w_j^n\in c_j^n$ such that $d_n(w_1^n,w_2^n)<2\ep$.
Furthermore, the diameters of the components of $\g_n\setminus\{w_1,w_2\}$ are bounded from below
by a constant which depends only on $r$ and $\diam\, \g_n$. This ends the proofs
of the claim and of the theorem.\qed

\subsection{Topological mixing properties}

We will assume throughout this subsection that $f$ is a
topologically cxc map
from the sphere to itself, 
endowed with a metric $d_v$ as in Theorem \ref{thm:can_gauge}. 
The main result, Proposition \ref{coro:univ}, 
articulates in quantitative form the principle of the Conformal Elevator: 
arbitrarily small balls can be blown up via the dynamics with controlled distortion.
The proof will use the fact (\cite[Cor.\,3.5.4]{kmp:ph:cxci}) that the map $f$ is a so-called
metric cxc map with respect to $(S^2,d_v)$.
Before recalling the properties needed, we define  the notion of roundness. 

\gap

\noindent{\bf Roundness.} Let $(Z,d)$ be a metric space and let $A$ be a bounded, proper  subset of $Z$ with nonempty interior.  Given $a \in \interior(A)$, let 
\[ L(A,a)=\sup\{d(a,b): b \in A\}\]
and
\[ l(A,a)=\sup\{ r : r \leq L(A,a) \; \mbox{ and } \; B(a,r) \subset A\}\]
denote, respectively, the {\em outradius} and {\em inradius} of $A$ about $a$.  
While the outradius is intrinsic, the inradius depends on how $A$ sits in $Z$. 
The condition $r \leq L(A,a)$ is necessary to guarantee that the outradius is at least the inradius. 
The {\em roundness of $A$ about $a$} is defined as 
\[ \roundness(A,a) = L(A,a)/l(A,a) \in [1, \infty).\]
One says $A$ is {\em $K$-almost-round} if $\roundness(A,a) \leq K$ for some $a \in A$, 
and this implies that for some $s>0$, 
\[ B(a,s) \subset A \subset B(a,Ks).\]
\gap

Being metric cxc implies the existence of 
\bi
\item continuous, increasing embeddings $\rho_{\pm}:[1,\infty) \to [1,\infty)$ such that,
for all $n, k \in \N$ and for all
$\wtU \in \UUU_{n+k}$, $\tilde{y} \in \wtU$, if $ U = f^{k}(\wtU)\in \UUU_n$ and $y=f^{k}(\tilde{y})$
then 
$$ \roundness(\wtU, \tilde{y}) <\rho_-(\roundness(U,y)) \quad \hbox{and} \quad
\roundness(U, y) <
\rho_+(\roundness(\wtU,
\tilde{y}))\,;$$
\item increasing homeomorphisms $\delta_{\pm}:[0,1] \to [0,1]$ such that, for all $n_0, n_1, k \in \N$ and for all
$\wtU \in \UUU_{n_0+k}$, $\wtU'\in\UUU_{n_1+k}$ with $\wtU' \subset \wtU$, if
$U =f^k(\wtU)\in \UUU_{n_0}$ and $U'=,f^k(\wtU') \in \UUU_{n_1}$, 
then
\[ \delta_+^{-1}\left(\frac{\diam \,U'}{\diam\,U}\right)\le \frac{\diam\wtU'}{\diam\wtU} 
\le \delta_-\left(\frac{\diam U'}{\diam
U}\right)\,.\]
\ib

\gap

To control distortion, some ``Koebe space'' is needed.  
To this end, let $\PPP_{k,K}$ denote the set of {\em preferred pairs}  $(W',W)$ of elements of {\bf U} such that
$|W'|=|W|+k$, $W'\subset W$ and $\roundness (W,x)\le K$ for any $x\in W'$.

\begin{prop}\label{prop:mix}
\be
\item For any $K\ge 1$, there exists $K'$ such that, whenever $(U',U)\in\PPP_{k,K}$, $W',W\in {\bf U}$ are
such that $f^p(U')=f^q(W')$ and $f^p(U)=f^q(W)$ belong to {\bf U} for some iterates $p$ and $q$, and $W'\subset W$,
then $(W',W)\in\PPP_{k,K'}$.
\item There exist $n_1$ and $K$ such that, for any $k_j\in\N$, any $U_j\in \UUU_{k_j}$, $j=1,2$, there exists $\wtU_1 \in f^{-(k_2+n_1)}(\{U_1\})$
such that $(\wtU_1,U_2)\in \PPP_{k_1+n_1,K}$.
\item There exists $k_0$ and $K_0$ such that, for any $n\ge 0$ and any $k\ge k_0$,
\be
\item for any $W'\in\UUU_{n+k}$, there exists $W\in\UUU_n$ such that $(W',W)\in\PPP_{k,K_0}$;
\item for any $W\in\UUU_n$, there exists $W'\in\UUU_{n+k}$ such that $(W',W)\in\PPP_{k,K_0}$.
\eb 
\eb
\end{prop}

\pf The first statement follows from the roundness distortion bounds enjoyed by metric cxc maps, by choosing $K'=(\rho_-\circ\rho_+)(K)$.
The conclusion in case (a) of the last statement follows directly from the diameter and roundness bounds using the fact that $\UUU$ is finite;
case (b) follows from the second point.

We now prove the second point, and start by fixing some notation.
There exists $r_i>0$ such that, for any $U\in\UUU$, there is some $x_U\in U$ with $B(x_U,r_i)\subset U$.
Let $n_0\ge 1$ be such that $f^{n_0}(U\cap X)= X$ for any $U\in\UUU$.
Choose $m\ge n_0$ as small as possible such that $2d_m\le r_i/2$.

We first assume that $U_2\in \UUU$. Let $W\in\UUU_{m}$ contain $x_{U_2}$. Since $f^{2m}(W)= X$, there is a component $\wtU_1$ of $f^{-2m}(U_1)$ which intersects $W$. It follows that $B(y,r_i - 2d_m)\subset U_2$ for all $y\in\wtU_1$ since, for all $z\in B(y,r_i-2d_m)$,
\begin{eqnarray*} |x_{U_2}-z|& \le &   |x_{U_2}-y| + |y-z|\\
& <  & \diam\,W +\diam\,\wtU_1 + (r_i-2 d_m)\\
&  \le & r_i.\end{eqnarray*}
Therefore, $\roundness(U_2,y)\le 2d_0/r_i$.
Thus, $(\wtU_1,U_2)\in\PPP_{2m + k_1, 2d_0/r_i}$.

Pick now $U_2$ randomly. It follows from above that there exists $\wtU_1'\in\UUU_{2m+k_1}$ such that $(\wtU_1',f^{k_2}(U_2))\in\PPP_{2m + k_1, 2d_0/r_i}$. By choosing a connected component $\wtU_1$ of $f^{-k_2}(\wtU_1)$ in $U_2$, we obtain $(\wtU_1,U_2)\in\PPP_{n_1+k_1,K}$
with $n_1=2m$ and $K=\rho_-(2d_0/r_i)$.
\qed

We derive the following property:

\begin{cor}[Injective conformal elevator]\label{coro:univ}
For any $X\in\GGG(f)$, there exist a distortion function $\eta_{ice}$, constants $c >0$ and $r_0>0$ such that,
for any $x\in X$ and $r >0$, there are an iterate $n\ge 0$ and a ball $B\subset B(x,r)$
of radius at least $c\cdot r$ such that $f^n|_{B}$ is $\eta_{ice}$-quasisymmetric and $f^n(B)$ contains
a ball of radius at least $r_0$. \end{cor}

\pf It is enough to prove the statement for the metric $d_v$.
Recall that $\mathbf{U}$ denotes the countable set of components of preimages of elements of $\UUU_0$ under iterates of $f$. 
The degree hypothesis in the definition of topologically cxc implies that we can find some 
$W \in \UUU_{k}$ such that the degree of $f^{k}:W\to f^k(W)$ is maximal, 
so that any further preimages $\wtW$ of $W$ map onto $W$ by degree one i.e., are homeomorphisms. 

It follows from \cite[Prop.\,3.3.2]{kmp:ph:cxci} 
that $W$ contains some ball $B(\xi,4 R)$ such that  for any iterate $n$, any $\txi\in f^{-n}(\xi)$,
the restriction $f^n:B(\txi,4R\te^{ n})\to B(\xi,4R)$ is a homeomorphism. 
\cite[Prop.\,3.2.3]{kmp:ph:cxci} shows that 
$f^n:B(\txi,R\te^{ n})\to B(\xi,R)$ is a similarity (in particular, it is quasisymmetric) and that the diameter of the unique preimage $\wtW$ of $W$ containing $B(\txi, R\te^n)$ is comparable to $R\te^n$.

Now suppose $B(x,r)$ is an arbitrary ball.  By \cite[Prop.\,2.6.6]{kmp:ph:cxci}, there exists $U \in \mathbf{U}$ with $U\subset B(x,r)$ and $\diam\, U\gtrsim r$.
But Proposition \ref{prop:mix} implies the existence of an iterated  preimage $\wtW$ of $W$ with $\wtW\subset U$ and $||\wtW| - |U||= O(1)$.  So $U$ contains $\wtW$, which is not too small; the previous paragraph implies $\wtW$ contains the desired ball $B$.
\qed

\subsection{Density points}

Let $(Z,d,\mu)$ be a $Q$-regular metric space. 
We recall that $\mu$ extends to an outer measure (which we also denote by $\mu$) on the power set of $Z$ so that if $A$ is any  subset of $Z$,
then there exists a Borel set $A^*$ containing
$A$ such that $\mu(A)=\mu(A^*)$.

A point $a\in A$ is an {\em $m$-density point} if $$\lim_{r\to 0}\frac{\mu(A\cap B(a,r))}{\mu(B(a,r))}=1\,.$$
A point $a\in A$ is a {\em $t$-density point} if, for all $\varepsilon >0$,  $$\lim_{r\to 0}\frac{1}{r}\sup_{z\in B(a,r)} d(z,A\cap B(a,(1+\varepsilon)r))=0\,.$$
Informally: the set $A$ becomes hairier and hairier upon zooming in at a $t$-density point. 

\begin{lemma}\label{lma:mt_density}  Suppose $Z$ is doubling and $A \subset Z$.  Then each $m$-density point of $A$ is also a $t$-density point.\end{lemma}

\pf If not, there are an $m$-density point $a$ of $A$, positive constants $\varepsilon$ and $c>0$, sequences of radii $(r_n)$ tending to $0$ and of points $z_n\in B(a,r_n)$
such that $d(z_n,A\cap B(a,(1+\varepsilon)r_n))\ge c r_n$. Of course $c\le 1$ since $a\in A$, so $B(z_n,\varepsilon c r_n)\subset B(a,(1+\varepsilon)r_n)\setminus A$ and
$$\frac{\mu ( B(a,(1+\varepsilon)r_n)\setminus A)}{\mu(B(a,(1+\varepsilon)r_n))}\ge \frac{\mu ( B(z_n,\varepsilon c r_n))}{\mu(B(a,(1+\varepsilon)r_n))} \gtrsim 1$$
which contradicts that $a$  is an $m$-density point of $A$.\qed

If $A \subset Z$ and $a\in A$, any tangent space to $(A, d|_A)$ at $a$ is naturally a subspace of a tangent 
space to $(Z, d)$ at $a$.  Equality holds if $a$ is a $t$-density point of $A$.

\begin{lemma}\label{lma:density_tangent}
Let $a$ be a $t$-density point of a subset $A$ of a doubling proper metric space $Z$. 
Suppose $r_n \downarrow 0$ is a decreasing sequence of radii such that the spaces $(Z_n, d_n=d/r_n, a)$ 
converge to a tangent $(T, d, a)$. 
Let $A'\subset T$ be the set of points $z$ such that there are $a_n\in Z_n\cap A$ such
that $a_n \to z$. Then $A'=T$.\end{lemma}

\pf Let $z\in T$, and consider a sequence of points $z_n$ which tends to $z$.
Since $a$ is a $t$-density point, for each $n$ we have 
$$d(z_n, A) \le  \sup_{x\in B(a,|z_n-a|)} d(x,A\cap B(a,2|z_n-a|))$$
so there is a point $a_n$ in $A\cap B(a,2|z_n-a|)$ such that $$\lim_n\frac{|z_n-a_n|}{|z_n-a|}=0\,.$$
Therefore, $$ d_n(z_n,a_n)\le \frac{|z_n-a|}{r_n}\cdot \frac{|z_n-a_n|}{|z_n-a|}\,.$$
By definition, the first fraction is uniformly bounded and the second tends to $0$: this 
proves that $z\in A'$.\qed

\subsection{Weak tangent spaces for topologically cxc maps}

We now analyze weak tangents to metric spaces in the conformal gauge of a
topologically cxc map $f$.  
Let $d_v$ be a visual metric given by Theorem \ref{thm:can_gauge}. 
The metric space $(S^2, d_v)$ is doubling, a property preserved by quasisymmetries.  
Hence for any $d\in\GGG(f)$, any sequence $r_n\to 0$, and any sequence $(x_n)$
of points in $X=(S^2, d)$, one may find a subsequence such that
$(X, d/r_n,x_n)$ tends in the
Gromov-Hausdorff topology to a doubling metric space $(T,t)$.

\begin{prop}\label{prop:tgt} For any $X\in \GGG(f)$, and any weak tangent
space $(T,t)$ of $X$, there is an open, onto
map $h:T\to X$ with discrete fibers, and  there is some constant $0 < c < 1$ such that, 
for any $R>0$, there is some $r_0$ such that any ball 
$B(x,r)$ with $x\in B(t,R)$ and $r\in(0,r_0)$ contains a ball
$B$ of radius $cr$ such that $g|_B$ is $\eta$-quasisymmetric. 
\end{prop}

We will use the following lemma whose proof is left to the reader.

\begin{lemma}\label{lma:tgt_equiv} Let $X\in\GGG(f)$. Then any tangent space of
$X$ is quasisymmetrically equivalent to a tangent
space of a visual metric metric $(S^2, d_v)$. \end{lemma}

\pf  (Prop. \ref{prop:tgt}) Lemma \ref{lma:tgt_equiv} implies
that it is enough to treat the case $d=d_v$.

By assumption $(T,t)$ is a limit of $(X_n,x_n)$ where $X_n=(X, d_v/r_n)$
and $(r_n)$ tends to zero.

Axiom [Irred] guarantees an integer $n_0$ so $f^{n_0}(U)=X$ for all $U \in \UUU_0$.  
By Proposition \ref{prop:mix}(3), there exists $K>1$ such that for each $n$ sufficiently large, 
there exists $W_n\in {\bf U}$ containing $x_n$ such that
$\roundness(W_n,x_n)\le K$ and $L(W_n,x_n)\asymp r_n$; that is, $W_n$ is $K$-almost round about $x_n$ and 
has diameter comparable to $1$ in the rescaled metric $d_v/r_n$.
Let $k_n=|W_n|$.
>From the diameter estimates (3) of Theorem \ref{thm:can_gauge}, the family of maps 
$$(f^{k_n+n_0}: X_n\to X)_n$$
is uniformly Lipschitz for some constant independent of $n$. 
Passing to a subsequence, we may assume the sequence $(f^{k_n+n_0}: X_n\to X)_n$ converges 
(in the Gromov-Hausdorff sense)
to a Lipschitz map $h:T\to X$. Since 
$f^{k_n+n_0}(W_n)=X$ by construction, the limit $h$ is also onto.

Let us prove that $h$ is open. Pick $z\in T$ and $r \in (0,\diam X)$, and choose
$z_n\in X_n$ with $z_n\to z$. By the diameter estimates (3) of Theorem \ref{thm:can_gauge}, 
$f^{k_n+n_0}(B_{X_n}(z_n,r))$ is ball centered at
$f^{k_n+n_0}(z_n)$ with definite radius, at least some uniform constant $r'>0$. It follows that, for $n$ large enough,
it covers $\overline{B(h(z),r'/2)}$, so $h(B(z,r))\supset B(h(z),r'/2)$.

Similarly, Corollary \ref{coro:univ} implies 
the quasisymmetry property; we leave the details to the reader.

Lastly, we prove that $h$ has discrete fibers. Fix a large ball $B(t,R) \subset T$. Pick $W^R_n\in  {\bf U}$ 
containing $x_n$ such that
$\roundness(W^R_n,x_n)\le K$, $B(x_n,R\cdot r_n)\subset W^R_n$ and $L(W_n,x_n)\asymp R\cdot r_n$.
For each $n$, $||W^R_n|-|W_n||= O(\log R)$ holds so that $f^{k_n+n_0}|_{W^R_n}$
have multiplicity   bounded by a constant depending only on $R$. This implies that $h|_{B(t,R)}$
has also bounded multiplicity.

\qed

Applying a theorem of Whyburn  \cite[Thm.\,X.5.1]{whyburn:analytic_topology}, Theorem \ref{thm:limit_surface} implies

\begin{cor}
\label{cor:tangent_map_surface}
If $f: S^2 \to S^2$ is topologically cxc, and $X \in \GGG(f)$,
then the map $h: T \to X$ given by Proposition \ref{prop:tgt} is a branched covering: 
it is locally of the form $z \mapsto z^k$ for some integer $k \geq 1$ and the branch set is discrete.  
\end{cor}

\section{Moduli of curves}

This section discusses curves and different notions of moduli
of family of curves.

\subsection{Rectifiable curves and curve families}
\label{subsecn:rectifiable}

Let $X$ be a proper metric space.  
In this subsection, we discuss some generalities on curves, see
e.g \cite[Chap.\,1]{vaisala:lectures_qc} for details.  

A {\em curve} is a continuous map  $\g:I \to X$ of a possibly infinite interval $I$ 
into $X$; 
it is {\em degenerate} if it is constant, and {\it nondegenerate} otherwise. 
Given two curves $\g:I\to X$ and $\g':I'\to X$,
we say they {\it differ by a reparametrization} if there exists a (not necessarily strict) monotone 
and onto map $\al:I\to I'$ such that $\g=\g'\circ\al$ or $\al:I'\to I$ such that $\g'=\g\circ\al$.
Note that two curves which differ by a reparametrization are simultaneously rectifiable or not.

We want to rescale curves of possibly variable length; for this, it is useful to think of such curves 
as maps of a common domain $\R$ into $X$.
Let $\g:[0,\ell(\g)]\to X$ be a rectifiable curve parametrized by arclength.
Given $t_0\in [0,\ell(\g)]$ and $x_0=\g(t_0)$, we define its {\it complete parametrization at $t_0$}, 
denoted  $(\g^0,t_0)$
as the map $\g^0:(\R, 0)\to (X, x_0)$ defined by
$\g^0(t)=\g(t+t_0)$ on $[-t_0,\ell(\g)-t_0]$, $\g_0|_{(-\infty,-t_0]}=\g(0)$ and 
$\g^0|_{[\ell(\g)-t_0, \infty)}=\g(\ell(\g))$.
We let $I_{\g^0}$ be the minimal compact interval such that $\g^0(I_{\g^0})=\g^0(\R)$.  
Note that given $x\in \g(I)$, there might be several complete parametrizations $\g^0$ such that
$\g^0(0)=x$.

A sequence $(\g^0_n: (\R, 0) \to (X, x_n))_n$ of completely parametrized curves {\em converges} 
if it converges uniformly on compact subsets of $\R$; note that this implies that the sequence $x_n$ converges. 
The limit of such a sequence is a locally rectifiable curve.
Given a bounded subset $A$ of $X$, the set of completely parametrized curves $(\R, 0) \to (X, x)$ for which $x \in A$ 
is relatively compact.

Given a completely parametrized curve $\gamma:(\R, 0) \to (X, x_0)$ and a small constant $r>0$, 
the rescaled curve given by $t\mapsto \g(t/r)$ defines again a completely parametrized curve 
$\gamma/r: (\R, 0) \to ((1/r)X, x_0)$.  Note that the intervals on which these rescaled curves are nonconstant 
satisfy $I_{\g/r}=(1/r)I_{\g^0}$ and therefore grow to $\R$ as $r \to 0$.
\gap

The following proposition comes from  \cite[Lma\,9.1]{margulis:mostow}.

\begin{prop}\label{prop:rect} Let $\g:I\to X$ be a rectifiable curve parametrised by arc length
in a metric space $X$. Then
$$\lim_{h\to 0} \frac{d(\g(t+h),\g(t-h))}{2|h|}=1$$
for almost all $t\in I$.
\end{prop}

Note that if $$\lim_{h\to 0} \frac{d(\g(t+h),\g(t-h))}{2|h|}=1$$
holds for some $t \in I$ then we also have
$$\lim_{h\to 0} \frac{d(\g(t+h),\g(t))}{|h|}=1$$
since $$\frac{d(\g(t+h),\g(t-h))}{2|h|}\le \frac{1}{2}\left( \frac{d(\g(t+h),\g(t))}{|h|}+\frac{d(\g(t-h),\g(t))}{|h|}\right)\le 1$$
and each term of the sum is also bounded by $1$.  That is, rectifiable curves are asymptotically  geodesic near almost every point. 

The next proposition illustrates how this property of rectifiable curves will be used
to produce many geodesics on tangent spaces. The hypothesis gives some uniformity to the rate at which 
rectifiable curves become asymptotically geodesic. 

\begin{prop}\label{coro:gdq} Let $(X,d)$ be a doubling proper metric space and $(\ep_k)$ be a sequence
of positive reals tending to $0$.  Suppose $A \subset X$ and for each $x \in A$, 
 there exists a completely parametrized curve
$\g: (\R, 0) \to (X, x)$ such that, for all $k$,  $$\sup_{|h|\le 1/k}\left|\frac{d(\g(h),\g(-h))}{2|h|}-1\right|\le \ep_k\,.$$

Assume that $x_0$ is a point of $t$-density of $A$.
Let $(r_n)$ be a sequence of positive reals tending to zero, and
assume that the limit of the rescaled spaces $X_n=(X,x_0, (1/r_n) d)$
exists in the pointed Gromov-Hausdorff topology. Let $(Z,z_0)$ be the limit.

Then any point in $Z$ lies in the image of a bi-infinite geodesic curve
$\gamma: \R \to Z$ which is a limit of rescaled completely parametrized curves passing through points of $A$.
\end{prop}

We say that the  geodesic curve $\gamma$ was {\it obtained by blowing up}.
The rescaled curves will be completely parametrized curves in $X_n$
which are rescalings of completely parametrized curves in $X$.

\gap

\pf Let $z\in Z$; by Lemma \ref{lma:density_tangent},
since $x_0$ is a $t$-density point,
one may find a sequence of points
$x_n\in X_n\cap A$ which tends to $z$. Consider complete parametrized curves
$(\hat{\g}_n,x_n)$ in $X$. Let $\g_n(t)=\hat{\g}_n(r_n t)$ be its complete parametrization
in $X_n$.  Since the original parametrizations are by arc length,
each $\g_n$ is a rectifiable curve which is represented by a $1$-Lipschitz function:
let $t,t'\in\R$,
$$ (1/r_n) d(\g_n(t),\g_n(t'))=(1/r_n) d(\hat{\g}_n(t r_n),\hat{\g}_n(t' r_n)) \le (1/r_n )|(t r_n) -(t'r_n)|\,.$$
Thus, Arzela-Ascoli's theorem implies the existence of a $1$-Lipschitz limit
$\g:\R\to Z$, with $\g(0)=z$.

We shall prove that $\g$ is a geodesic.

By the choice of $\gamma_n$ we have
$$\lim_{h\to 0} \sup_n\left|\frac{d(\hat{\g}_n(h),\hat{\g}_n(0))}{|h|}-1\right|=0\,,$$
and
$$\lim_{n\to\infty} (1/r_n ) d(\g_n(t),\g_n(0))= \lim_{n\to\infty} |t|\frac{d(\hat{\g}_n(t r_n),\hat{\g}_n(0))}{|t| r_n}= |t|\,.$$
Hence $d(\g(t),\g(0))=|t|$.
Similarly,
$$\lim_{n\to\infty} (1/r_n ) d(\g_n(t),\g_n(-t))= \lim_{n\to\infty} |t|\frac{d(\hat{\g}_n(t r_n),\hat{\g}_n(-t r_n))}{|t| r_n}= 2|t|$$
and $d(\g(t),\g(-t))=2|t|$.

Let $t,t'$ be non-zero real numbers. We may assume that $|t|\ge |t'|$.

We first assume that they have the same sign. Hence $|t-t'|= |t|-|t'|$. Therefore,
$$|t-t'|= | d(\g(t),\g(0))- d(\g(t'),\g(0))|\le d(\g(t),\g(t'))\le |t-t'|\,.$$

If $t$ and $t'$ have opposite signs, then $|t-t'|= |t|+|t'|$ and $|t+t'|=|t|-|t'|$.
Thus,
$$|t-t'|=2|t|-(|t|-|t'|)=d(\g(t),\g(-t))- d(\g(t'),\g(-t))\le d(\g(t),\g(t'))\le |t-t'|\,.$$

Therefore, $\g$ is geodesic.\qed

\subsection{Analytic moduli} \label{subscn:analyticmoduli}

Suppose $(X,d,\mu)$ is a metric measure space,
$\Gamma$ is a family of curves in $X$, and $p\ge 1$.  The {\em (analytic) $p$-modulus} of $\Gamma$ is defined by 
$$\mod_p\Gamma=\inf \int_X \rho^p d\mu$$
where the infimum is taken over all Borel functions
$\rho:X\to[0, +\infty] $ such that $\rho$ is {\it admissible}, i.e. $$\int_{\gamma}\rho \ ds\ge 1$$ for all
$\gamma\in \Gamma$ which are rectifiable. 
If $\Gamma$ contains no rectifiable curves, $\mod_p\Gamma$ is defined to be zero.  
Note that when $\Gamma$ contains  a constant curve, then
there are no admissible $\rho$, so we set $\mod_p\Gamma=+\infty$.
We say that a family of curves is {\it nondegenerate} if it contains no constant curves.
Modulus behaves like an outer measure: it is countably subadditive, etc.---see \cite{heinonen:koskela:loewner}.

We note that the moduli of two families of curves $\G$ and $\G'$ are the same if each curve of $\G'$ differs from
a curve of $\G$ by reparametrization, and vice-versa.

\gap

We recall the following basic estimate in
an Ahlfors-regular metric space \cite[Lma\,3.14]{heinonen:koskela:loewner}:

\begin{prop}\label{prop:upperbound} Let $X$ be a $Q$-Ahlfors regular
metric space with $Q>1$, $\G$ be the family of curves which joins
$B(x,r)$ to $X\setminus B(x,R)$ for some $x\in X$, and $0<r<R\le \diam\,X$.
Then $$\mod_Q\G\lesssim \frac{1}{\log^{Q-1} (R/r)}\,.$$
In particular, 
the family of non-trivial curves 
which go through $x$
has zero $Q$-modulus.\end{prop}

\subsection{Thick curves}

Suppose $X$ is a $Q$-Ahlfors regular metric space.  
Denote by  $\PPP(X)$ the set of compact continuous curves $\gamma: [0,1] \to X$, 
endowed with the supremum-norm topology inherited from the metric of $X$.  
Thus, if $\gamma$ is a curve, then $B_\infty(\gamma,r)$ denotes the set of curves 
$\g':[0,1]\to X$
such that $|\g(t)-\g'(t)| < r$ for all $t\in [0,1]$. 

\gap

Following Bonk and Kleiner, we say a curve $\gamma \in \PPP(X)$ is {\em thick} if it is nonconstant and the $Q$-modulus
of $B_\infty(\gamma,\epsilon)$ is positive for all $\epsilon >0$ \cite{bonk:kleiner:conf_dim}. 
The following properties are established by Bonk and Kleiner:

\begin{prop}\label{prop:thick}
Let $X$ be an Ahlfors regular space. The set of non-thick curves has zero modulus. 
A nonconstant limit in $\PPP(X)$ of thick curves is thick, and nonconstant connected subcurves of thick curves are
thick.
\end{prop}

If there exist curve families of positive modulus, the preceding proposition implies the existence of many thick curves.  The next proposition refines and localizes this result.

\begin{prop}\label{prop:thick2} Let $X$ be an Ahlfors regular space and $\G\subset\PPP(X)$ a family
of nondegenerate curves of positive modulus. The subfamily $\G'$ of curves
$\g\in\G$ such that there is some $\ep>0$ with the property that $\mod_Q(B_\infty(\g,\ep)\cap\G)=0$
has zero modulus.\end{prop}

\pf We may cover $\G'\subset\PPP(X)$ by balls $B_i$ of uniformly bounded radius such that
$\mod_Q B_i\cap \G=0$. Since $\PPP(X)$ is separable, we may extract a countable subcover:
the proposition follows from the $\si$-subadditivity of modulus.\qed

\gap

Let us observe the following fact: let $\g:I\to X$ and $\g':I'\to X$ be two curves which differ
by a reparametrization. Given $\ep>0$, denote by $\G_\ep$ (resp. $\G_\ep'$) the set of curves $c:I\to X$ 
(resp. $c':I'\to X$) such that $\sup_{t\in I}|c(t)-\g(t)|\le \ep$ (resp. $\sup_{t\in I'}|c'(t)-\g'(t)|\le \ep$).
Let $\al:I\to I'$ be a monotone onto map such that $\g=\g'\circ\al$. If $c'\in\G_\ep'$, then $(c'\circ\al)\in\G_\ep$.
Conversely, assume that $c\in\G_\ep$. Since $\al$ is monotone, it is a uniform limit of homeomorphisms
$\al_n:I\to I'$. Then, for all $t\in I'$, 
$$|\al\circ\al_n^{-1}(t) - t|=|\al\circ\al_n^{-1}(t) - \al_n\circ\al_n^{-1}(t)|\le
\|\al-\al_n\|_\infty\,.$$ Let $\om$ be a modulus of continuity for $\g'$. 
It follows that 
\begin{eqnarray*} |c\circ\al_n^{-1}(t)-\g'(t)|& \le & |c\circ\al_n^{-1}(t)-\g\circ\al_n^{-1}(t)| + |\g\circ\al_n^{-1}(t)-\g'(t)|\\
&\le & \|c-\g\|_\infty + |\g'(\al\circ\al_n^{-1}(t))- \g'(t)|\le \ep +\om (\|\al-\al_n\|_\infty)\,.
\end{eqnarray*}
Therefore, if $n$ is large enough then $c\circ\al_n^{-1}\in\G_{2\ep}'$. 
>From this discussion, we see that thickness is a notion which may be generalized to curves defined on other intervals
than $[0,1]$, and that two curves which differ by a reparametrization are simultaneously thick or not.

\gap

We record here the following result:

\begin{prop}[Tyson]\label{prop:tyson} Let $X$ and $Y$ be two regular metric
spaces of dimension $Q>1$ and $h:X\to Y$ be a quasisymmetric map. 
There exists a constant $K\ge 1$ such that, for any
family of curves $\G\subset X$, $$\frac{1}{K}\mod_Q\G\le \mod_Q h(\G)\le K\mod_Q\G\,.$$
In particular, quasisymmetric maps preserve the set of thick curves.\end{prop}

For a proof, see \cite[Thm.\,1.4]{tyson:qcqs} (see also Prop.\,\ref{prop:fthick=thick}).

\gap

Note that thick curves need not exist: if there is no family of curves of positive modulus,
which happens for instance if there are no rectifiable curves, then there are no thick curves.

\subsection{Support of curve families}

If $\G$ is family of curves, we define its support $\supp\,\G$ as the set of points
$x\in X$ such that $x\in \g$ for some $\g\in\G$.
We note that $\supp\,\G$ might not be measurable.
In any case, $\supp\,\G$ is contained in some
Borel set of the same measure. We will then denote by $\supp^*\G$ such a set.

\gap

For $L>0$, 
let $\Gamma_L = \{\gamma \in \Gamma,\  \ell(\gamma) \leq L\}$ be the subfamily of curves in $\Gamma$ of 
length at most $L$; note that every element of $\Gamma_L$ is rectifiable.  We say that $\Gamma$ is {\em closed} if, 
for all $L > 0$, the set of all complete parametrizations of elements of $\Gamma_L$ is closed in the supremum norm 
on functions $\R \to X$.  That is, whenever a sequence $(\gamma_n) \subset \Gamma_L$ has the property 
that some choice of complete parametrizations $(\gamma_n^0)$ converges, then the limit is, up to reparametrization, 
a curve in $\G_L$.  Note that according to this definition the property of being closed depends only on the subfamily of rectifiable curves.

\begin{lemma}\label{lma:cpctsup} Let $X$ be a $Q$-regular metric space for which there exists a family of 
nondegenerate curves
of positive $Q$-modulus. Then there exists $L>1$ and a  closed family $\G$ of
rectifiable curves of diameter at least $1/L$ and of length at most $L$ and with compact support and positive 
modulus.\end{lemma}

\pf According to Proposition \ref{prop:thick}, $X$ contains a rectifiable thick curve $\g \in \mathcal{P}(X)$.
Choose $r$ so that $0<r< \diam\,\g/3$ and put $\G_0=B_\infty(\gamma, r) \subset \mathcal{P}(X)$. 
The family $\G_0$ has no trivial curves and each
curve has diameter at least $\diam\,\g/3$.
Define $\G_n$ as the curves of $\G_0$ of length at most $n$.
By the $\si$-subadditivity of moduli,
there is some $n$ such that $\mod_Q\G_n >0$.
The closure $\Gamma$ of the set of all complete parametrizations of elements of $\G_{n}$
in the sup-norm gives the desired family; the additional curves added as limit points are rectifiable.  
\qed

\subsection{Combinatorial moduli}

After recalling the definition of combinatorial moduli,
we establish and recall some basic estimates which will be used later on.

\subsubsection{Definitions and properties}

\noindent{\bf Definitions.}
\gap
Let $\SSS$ be a covering of a  topological
space $X$, and let $p\ge 1$.
Denote by $\MMM_p(\SSS)$ the set
of functions $\rho:\SSS\to \R_+$ such that $0<\sum\rho(s)^p<\infty$; 
elements of $\MMM_p(\SSS)$ we call  {\it admissible metrics}.
For $K\subset X$ we denote by $\SSS(K)$ the set of elements of $\SSS$ which intersect $K$.  
The {\it $\rho$-length} of $K$ is by definition
$$\ell_\rho(K)=\sum_{s\in\SSS( K)} \rho(s)\,.$$
Define the {\em $\rho$-volume} by
$$V_{p}(\rho)=\sum_{s\in\SSS} \rho(s)^p\,.$$
If $\Gamma$ is a family of curves in $X$ and if $\rho\in\MMM_p(\SSS) 
$, we define
$$L_\rho(\Gamma,\SSS)=\inf_{\gamma\in\Gamma} \ell_\rho(\gamma),$$
$$\mod_Q(\Gamma, \rho, \SSS)=\frac{V_{p}(\rho)}{L_\rho(\Gamma,\SSS)^p},$$ 
and the {\em combinatorial modulus} by
$$\mod_p(\Gamma,\SSS) = \inf_{\rho\in\MMM_p(\SSS)}
\mod_p(\Gamma,\rho,\SSS).$$

Note that if $\SSS$ is a finite cover, then the modulus of a nonempty family of curves
is always finite and positive.

A metric $\rho$ for which $\mod_p(\Gamma,\rho,\SSS)=\mod_p(\Gamma, \SSS)$ will be called {\em optimal}.  
We will consider here only finite covers; in this case the proof of the existence of optimal metrics 
is a straightforward argument in linear algebra.  
The following result is the analog of the classical Beurling's criterion
which  characterises optimal metrics.

\begin{prop}\label{prop:metric_opt}
Let $\SSS$ be a finite cover of a space $X$,
$\Gamma$ a family of curves and $p>1$.
An admissible metric $\rho$ is optimal if and only if there is a non-empty finite subfamily $ 
\Gamma_0\subset\Gamma$ and non-negative scalars $\lambda_{\gamma}$, $ 
\gamma\in\Gamma_0$, such that \be
\item for all $\gamma\in\Gamma_0$, $\ell_\rho(\gamma) = L_\rho(\Gamma, 
\SSS)$\,;
\item for any $s\in\SSS$, $$p\rho(s)^{p-1}=\sum \lambda_{\gamma}$$
where the sum is taken over curves in $\Gamma_0$ which  
go through $s$.\eb
Moreover, an optimal metric is unique up to scale, and one has
$$\mod_p(\G,\SSS)=\frac{1}{p}\sum_{\g\in\G_0} \lambda_\g\,.$$\end{prop}

For a proof, see Proposition 2.1 and Lemma 2.2 in \cite{ph:emp}.
\gap

\begin{prop}\label{prop:subadditivity}
Let $\SSS$ be a locally finite cover of a topological space $X$ and $p \geq 1$.

\be
\item If $\Gamma_1 \subset \Gamma_2$ then $\mod_p(\Gamma_1,\SSS) \leq \mod_p(\Gamma_2, \SSS)$.
\item If every curve of $\Gamma_1$ contains a curve of $\Gamma_2$ then $\mod_p(\Gamma_1,\SSS)\le\mod_p(\Gamma_2,\SSS)$.
\item Let $\Gamma_1,\ldots,\Gamma_n$ be a set of curve families in $X$.  Then $$\mod_p(\cup\Gamma_j,\SSS)\le \sum\mod_p(\Gamma_j,\SSS)\,.$$
\eb
\end{prop}

The proof is the same as the standard one for classical moduli (see for instance \cite[Thm.\,6.2, Thm.\,6.7]{vaisala:lectures_qc}) and so is omitted.

\subsubsection{Dimension comparison}

Let $X$ be a topological space endowed with a finite cover $\SSS$. Consider a curve family $\G$.
For $p\ge 1$, when considering admissible metrics $\rho$ for the $p$-modulus, one may always assume that
$L_\rho(\G,\SSS)=1$ and $\rho\le 1$. Thus, if $1\le p\le q$, then $\mod_q(\G,\SSS)\le\mod_p(\G,\SSS)$.
The following proposition improves in some cases this estimate:

\begin{prop}\label{prop:comb_dim_estimate}
If $1\le p\le q$, then, for all $\ep>0$, the following holds:
$$\mod_q(\G,\SSS)\le \left(\ep^{q-p}+ \frac{1}{\ep^p}\sup_{s\in\SSS}\mod_q(\G(s),\SSS)\right)\mod_p(\G,\SSS)$$
where $\G(s)$ denotes the subfamily of curves of $\G$ which go through $s$.\end{prop}

\pf Let $\rho$ be the optimal metric for the $p$-modulus of $\G$ such that $L_\rho(\G, \SSS)=1$. 
For $\ep >0$, set $\EEE(\ep)=\{s\in\SSS,\ \rho(s)\ge \ep\}$. It follows from the Markov 
inequality that
$$card\,\EEE(\ep)\le \frac{1}{\ep^p}\mod_p(\G,\SSS)\,.$$
Decompose $\G$ into
two family of curves: those which go through at least one element of $\EEE(\ep)$, which we will denote
by $\G_{\EEE}$, and its complement which we denote by $\G'$.
Then, by the elementary properties of moduli, 

$$\left\{\begin{array}{l}
\mod_q(\G',\SSS)\le  \dis\sum_{s\in\SSS(\G')} \rho^q(s)\le \dis\left(\sup_{\SSS\setminus\EEE(\ep)}\rho \right)^{q-p} \dis\sum_{s\in\SSS(\G')} \rho^p(s)\le \ep^{q-p}\mod_p(\G,\SSS)\,,\\ \,\\
\mod_q(\G_{\EEE},\SSS)\le \dis\sum_{s\in\EEE(\ep)}\mod_q(\G(s),\SSS)\le \left(\dis\frac{1}{\ep^p}\mod_p(\G,\SSS) \right)\dis\sup_{s\in\SSS}\mod_q(\G(s),\SSS)\,.\end{array}\right.$$
Hence
$$\mod_q(\G,\SSS)\le \left(\ep^{q-p}+ \frac{1}{\ep^p}\sup_{s\in\SSS}\mod_q(\G(s),\SSS)\right)\mod_p(\G,\SSS)\,.$$\qed

\begin{cor}\label{cor:comb_dim}
Let $(\SSS_n)$ be a sequence of finite coverings and $q > p\ge 1$. If we assume 
there is a sequence $(\eta_n)$ of positive
numbers tending to zero such that
$$\sup_{s\in\SSS_n}\mod_q(\G(s),\SSS_n)\le \eta_n\,,$$
then
$$\lim_{n\to\infty} \frac{\mod_q(\G,\SSS_n)}{\mod_p(\G,\SSS_n)} =0\,.$$\end{cor}

This is yet further evidence that combinatorial modulus behaves like Hausdorff measure, under the appropriate assumptions.

\gap

\pf Pick $\ep_n= \eta_n^{1/(2p)}$ so that $(\ep_n)$ converges to $0$, and apply Proposition \ref{prop:comb_dim_estimate} for each $n$. It follows
that 
$$ \frac{\mod_q(\G,\SSS_n)}{\mod_p(\G,\SSS_n)} \le \ep_n^{q-p}+ \frac{1}{\ep_n^p}\eta_n
\le\ep_n^{q-p}+ \eta_n^{1/2}\,.$$
Therefore, since $q>p$, $$\lim_{n\to\infty} \frac{\mod_q(\G,\SSS_n)}{\mod_p(\G,\SSS_n)}=0\,.$$\qed

\subsubsection{Transformation rules}

\begin{prop}\label{prop:pvalentmod} Let $X$ and $X'$ be two connected, Hausdorff and locally compact topological spaces, and $\SSS$, $\SSS'$ respectively be coverings by compact connected subsets.
Let $f:X'\to X$ be an onto, proper and continuous map such that
\bi
\item for every $s'\in\SSS'$, $f(s')\in\SSS$ and for every $s\in\SSS$, the set of connected components of $f^{-1}(\{s\})$
is a subset of $\SSS'$;
\item there exists an integer $d\ge 1$ such that, for every $s\in\SSS$, the set $f^{-1}(\{s\})$
has at most $d$ connected components.
\ib
Let $\Gamma'\subset X'$ and $\Gamma\subset X$ be two family of curves.
Then
\be
\item if, for every $\g'\in\Gamma'$, $f(\g')$  contains a curve $\g\in\Gamma$, then
$$\mod_p(\Gamma',\SSS')\le d\cdot \mod_p(\Gamma,\SSS)\,;$$
\item if every curve in $\Gamma$ contains the image of a curve of $\Gamma'$, then
$$\mod_p(\Gamma',\SSS')\ge \frac{1}{d^p}\mod_p(\Gamma,\SSS)\,.$$
\eb\end{prop}

\pf
\be
\item
Let $\rho$ be an admissible metric for $\Gamma$.
Set $\rho'=\rho\circ f$. If $\g'\in\Gamma'$, then let $\g\in\Gamma$ be a subcurve of $f(\g')$.
One has $f(\SSS'(\g'))=\SSS(f(\g'))$ so
$$\ell_{\rho'}(\g') \ge\sum_{s\in\SSS (f(\g'))} \rho(s)\ge \sum_{s\in\SSS (\g)} \rho(s)\ge L(\Gamma,\rho)\,.$$
On the other hand, 
$V_p(\rho')\le d\cdot V_p(\rho)$  so that $$\mod_p(\Gamma',\SSS')\le d\mod_p(\Gamma,\SSS)\,.$$

\item
Let $\rho'$ be an admissible metric for $\Gamma'$.
Set $$\rho(s)=\left(\sum_{f(s')=s}\rho'(s')^p\right)^{1/p}\,.$$

It follows that 
$$\rho(s)\ge \max_{f(s')=s}\rho'(s')\ge\frac{1}{d}\sum_{f(s')=s}\rho'(s')\,.$$
Therefore, if $\g\in\Gamma$ and $\g'$ is a curve the image of which is contained in $\g$, then,
as $\SSS'(\g')\subset \SSS'( f^{-1}(\g))$,
$$\ell_\rho(\g)\ge \frac{1}{d}\ell_{\rho'}(\g')\ge \frac{1}{d} L(\Gamma',\rho')\,.$$
Moreover, $V_p(\rho) = V_p(\rho')$  so that 
$$\mod_p(\Gamma',\SSS')\ge \frac{1}{d^p}\mod_p(\Gamma,\SSS)\,.$$
\eb\qed

\subsubsection{Analytic versus combinatorial moduli}

Under suitable conditions, the combinatorial moduli obtained from a sequence $(\SSS_n)$ of coverings can be used to approximate
analytic moduli on metric measure spaces. 

\gap
The approximation result we use requires the sequence of coverings $(\SSS_n)$ to 
be a {\em uniform family of quasipackings}.  

\gap
\begin{defn}[Quasipacking] A {\em quasipacking} of a metric space is
a locally finite cover $\SSS$ such that there is some constant
$K\ge 1$ which satisfies the following property.
For any $s\in\SSS$, there are two balls $B(x_s, r_s)\subset s\subset B(x_s, K\cdot r_s)$ 
such that the family $\{B(x_s, r_s)\}_{s\in\SSS}$ consists of pairwise disjoint balls.  
A family $(\SSS_n)$ of quasipackings is called {\em uniform } if the mesh of $\SSS_n$ tends to zero as $n \to \infty$ and the constant $K$ defined above can be chosen independent of $n$. 
\end{defn}

Uniform quasipackings are preserved under quasisymmetric maps  quantitatively.

\gap

The next result says that under appropriate hypotheses, analytic and combinatorial moduli are comparable.  

\begin{prop}\label{prop:comb_mod} Suppose $Q > 1$, $X$ is an Ahlfors $Q$-regular compact  
metric space, and $(\SSS_n)$ is a sequence of uniform quasipackings. Let
$\Gamma$ be a nondegenerate closed family of curves in $X$. 
Then either 
\be
\item $\mod_Q\Gamma=0$ and $\lim_{n \to \infty}\mod_Q(\Gamma,\SSS_n)=0$, or 
\item $\mod_Q\Gamma>0$, and there exist constants $C \geq 1$ and $N \in \N$ such that  for any  $n>N$,
$$\frac{1}{C}\mod_Q(\Gamma,\SSS_n)\leq \mod_Q\Gamma\leq C\mod_Q(\Gamma,\SSS_n).$$
\eb
\end{prop}

See Proposition B.2 in \cite{ph:emp}.

Suppose now $f: S^2 \to S^2$ is topologically cxc, let $\SSS_n:=\UUU_n, n = 0, 1, 2, \ldots$ 
be the sequence of open covers as in the definition, and let $d$ be a metric in the conformal gauge of $f$.  
Then, according to \cite[Thm.\,4.1]{kmp:ph:cxcIII}, $\{\SSS_n\}$ is a uniform sequence of quasipackings, 
so we may apply Proposition \ref{prop:comb_mod} 
to estimate analytic moduli using combinatorial moduli.

\subsection{Intersection of curves}

\begin{prop}\label{prop:regular_intersection} Let $X$ be a $Q$-Ahlfors regular compact metric space, $Q > 2$.  
If $\G \subset \mathcal{P}(X)$ is a closed nondegenerate curve family of positive $Q$-modulus, then the family 
$$\Gamma^\perp = \{\gamma' \in \PPP(X) \mid  \forall\g\in\G,\ \gamma' \cap\g\ne\emptyset\}$$
is also closed and nondegenerate, and $\mod_Q \Gamma^\perp=0$.
\end{prop}

\gap

\pf That $\Gamma^\perp$ is closed follows immediately from the definitions. 
If $\Gamma^\perp$ contained a constant curve then every curve of $\G$ would go through a given point,
implying $\mod_Q \Gamma=0$ 
by Proposition \ref{prop:upperbound}; hence $\Gamma^\perp$ is nondegenerate.

Given $\de>0$, denote by $\G_{\delta}^\perp$ the family of curves $c\in\G^\perp$ of diameter at least
$\delta$.
The $\si$-subadditivity of the modulus
implies that it is enough to prove that $\mod_Q\G_{\delta}^\perp=0$ for all $\delta>0$.
Fix $\de>0$, and note that $\G_{\delta}^\perp$ is closed, 
so we may estimate its modulus using quasipackings,
cf. Proposition \ref{prop:comb_mod}.

Let $\RRR_n$ denote a maximal $2^{-n}$-separated set of $X$, then $\SSS_n=\{B(x,1/2^n), x\in\RRR_n\}$
defines a uniform sequence of quasipackings. 
It follows from the $Q$-regularity of $X$ that 
$$\lim_{n\to \infty} \sup_{s\in\SSS_n}\mod_Q(\G_{\delta}(s),\SSS_n)=0\,,$$
where $\G_{\delta}(s)$ denotes the set of curves of diameter at least $\delta$
which intersects $s$.
This is just the discrete version of Proposition \ref{prop:upperbound}: given $s_0=B(x_0,1/2^n)\in\SSS_n$, one may
estimate $\mod_Q(\G_{\delta}(s_0),\SSS_n)$ by setting $\rho_n(s)=0$ if $dist(x_0,s) > \delta/2$ or $s=s_0$,
and $\rho_n(s)= \diam\,s/dist(x_0,s)$ otherwise.

Fix $n\ge 1$.
By Proposition \ref{prop:metric_opt}, there are an optimal metric $\rho_n$ for 
$\mod_2(\Gamma_\de^\perp,\SSS_n)$, 
a nonempty subfamily $ \Gamma_n\subset\Gamma_\de^\perp$ and nonnegative scalars $\lambda_{\gamma}$, 
$\gamma\in\Gamma_n$, such that 
\be
\item for all $\gamma\in\Gamma_n$, $\ell_{\rho_n}(\gamma) = L_{\rho_n}(\Gamma_\de, 
\SSS_n)=1$\,;
\item for any $s\in\SSS_n$, $$2\rho_n(s)=\sum \lambda_{\gamma}$$
where the sum is taken over curves in $\Gamma_n$ which  
go through $s$.\eb

We note that any curve of $\Gamma$  intersects each curve of $\Gamma_n$ so that
$$L_{\rho_n}(\Gamma,\SSS_n)\ge \frac{1}{2}\sum_{\g\in\Gamma_n} \lambda_{\g}=v_2(\rho_n)\,.$$
It follows that $$\mod_2(\G,\SSS_n)\cdot \mod_2(\G_\de^\perp,\SSS_n)\le \left(\frac{v_2(\rho_n)}{v_2(\rho_n)^2}\right)v_2(\rho_n)\le 1\,.$$

According to Proposition \ref{prop:comb_mod}, one has, for any $n$ large enough,
$$0<\mod_Q\G\lesssim \mod_Q(\G,\SSS_n) \le \mod_2(\G,\SSS_n) $$ since $Q\ge 2$.
Therefore, $$\mod_2(\G_\de^\perp,\SSS_n)\lesssim \frac{1}{\mod_Q\G}<\infty$$ for all $n$.
Since $Q>2$, we may apply Corollary  \ref{cor:comb_dim} to $\G_\de^{\perp}$ with $q=Q$ and $p=2$ to conclude 
$$\lim_{n\to\infty}\mod_Q(\G_{\delta}^\perp,\SSS_n)=0\,.$$
\qed

Let $X$ be a metric space and suppose $\g, \g' \in \PPP(X)$ are two curves in $X$.  
We say that $\g$ and $\g'$ {\em cross} if there exists $\epsilon>0$ such that each curve of $B_\infty(\g, \epsilon)$ meets each curve of $B_\infty(\g', \epsilon)$.

\begin{cor}\label{cor:inter_surface_thick}
Let $X$ be a $Q$-regular metric space  with $Q\ge 2$. If there are two thick curves which
cross, then $Q=2$.\end{cor}

\pf Suppose $\g$ and $\g'$ are two thick curves.  
By definition, for all $\epsilon>0$, we have $\mod_Q B_\infty(\gamma, \epsilon)>0$ and  $\mod_Q B_\infty(\gamma', \epsilon)>0$.  
If $\g, \g'$ cross, then for some $\epsilon$ we have $B_\infty(\gamma', \epsilon) \subset B_\infty(\gamma, \epsilon)^\perp$.  
But Proposition \ref{prop:regular_intersection} implies $\mod_Q B_\infty(\gamma, \epsilon)^\perp =0$ and so 
$\mod_Q B_\infty(\gamma', \epsilon)=0$, which contradicts the thickness of $\gamma'$. 
\qed

We may now give a necessary and sufficient condition for a dynamical system
to be conjugate to a conformal one:

\begin{prop}\label{prop:mainthmrat}
Let $\DDD$ be either a topologically cxc mapping on $S^2$or a hyperbolic group with boundary $S^2$.
Assume that  the conformal dimension of $\DDD$ is attained by an Ahlfors $Q$-regular metric space $X$ in its conformal gauge. 
Then $\DDD$ is  topologically conjugate to a semihyperbolic rational map or cocompact Kleinian group if and only if there are
two thick curves in $X$ which cross.
\end{prop}

\pf 
Let $d$ be a metric of minimal dimension $Q\ge 2$. 

If $\DDD$ is conjugate to a genuine conformal dynamical system, then $Q=2$,
and there is a quasisymmetric map from $(X,d)$ to $\cbar$ by Theorem \ref{thm:can_gauge}.
Since the sphere admits thick curves which cross, so does $X$, according to Proposition \ref{prop:tyson}.

Conversely, if there exist thick curves in $X$ which cross, Corollary \ref{cor:inter_surface_thick} implies that $Q=2$.  
Theorem \ref{thm:char_of_rm} completes the proof. \qed


\section{Moduli at the minimal dimension}
In this section, we assume that $f: S^2 \to S^2$ is topologically cxc, and that $d$ is a $Q$-dimensional 
Ahlfors regular metric on $S^2$ belonging to the conformal gauge of $f$. We denote by $X=(S^2, d)$.  
Recall that the sequence $\SSS_n:=\UUU_n$ of coverings forms a uniform sequence of quasipackings.  

\subsection{Positive modulus on $X$}

\gap
\begin{prop} 
\label{prop:positive_modulus} 
There is a family of curves on $X$ with positive
$Q$-modulus.\end{prop}

\pf Since $Q$ is the AR-conformal dimension,
it follows from  \cite[Cor.\,1.0.2]{keith:laakso} 
that a weak tangent $(T,t)$ of $X$ admits a family $\Gamma$
of positive $Q$-modulus contained in some ball $B_T(t,R)$ for some $R>0$. 
We may assume that this family of curves is closed with compact support and definite diameter (Lemma \ref{lma:cpctsup}).
By Proposition \ref{prop:thick2}, we may also assume that $\mod_Q(B_\infty(\g,\ep)\cap\G)>0$
for every curve $\gamma$ in $\G$ and any $\ep>0$.
Let $c \in (0,1)$ be the constant given by Proposition \ref{prop:tgt}.  

Since $\G$ has positive modulus, its support $\supp\,\G$ has positive measure, so we may
find a point of $m$-density $z\in\supp\,\G$.  It follows that we may find a radius $r>0$ 
which is small enough so that any ball $B\subset B(z,r)$ of radius at least $cr/3$ satisfies 
$\mu(B\cap \supp\,\G) >0$. 

According to Proposition \ref{prop:tgt}, there exists a map $h:T\to X$ and a ball $B\subset B(z,r)$ of radius $cr$ 
such that $h:B\to X$ is a quasisymmetric embedding. Since $\mu(\supp\,\G\cap (1/3)B)>0$ by construction, there is some curve $\g\in \G$
which intersects $(1/3)B$ hence there is some $\ep>0$ such that every curve
from $\G\cap B_\infty(\g,\ep)$ intersects $(1/2)B$. It follows that $\mod_Q\G_0>0$
where $\G_0$ denotes the subcurves in $B$ of those curves in $\G$
which enter $(1/2)B$.

Therefore, from Proposition \ref{prop:tyson}, it follows that
$h(\Gamma_0)$ is a family of positive $Q$-modulus on $X$.\qed

\subsection{Invariance of thick curves}

\begin{prop}\label{prop:fthick=thick} The image of a thick curve under $f$ is a thick curve.\end{prop}

\pf Let $\g$ be a thick curve and let $\ep >0$. Since $f$ is uniformly continuous, there exists
$\delta>0$ such that $f(B_{\infty}(\g,\delta))\subset B_{\infty}(f(\g),\ep)$.
By definition $\mod_Q B_{\infty}(\g,\delta) >0$; there is some $L<\infty$ such that
the set $\G\subset B_{\infty}(\g,\delta)$ of curves of length at most $L$ has positive $Q$-modulus.
According to Proposition \ref{prop:comb_mod}, there exists $m>0$ such that
$$\mod_Q(\G,\SSS_n)\ge m >0$$ for all $n$ large enough. Therefore,
by Proposition \ref{prop:pvalentmod}, $\mod_Q (f(\G),\SSS_n)\gtrsim m$ holds, and another appeal to 
Proposition   \ref{prop:comb_mod} now implies that $$\mod_Q B_{\infty}(f(\g),\ep) \ge \mod_Qf(\G) >0\,.$$
\qed

The following corollary is needed for the construction of the second tangent space, $T_2$, mentioned in the introduction.  

\begin{cor}\label{cor:locallythick} If $T$ is a tangent space of $X$, 
$h:T\to X$ be given by Proposition \ref{prop:tgt},
and $\gamma: [0,1] \to T$ is a nonconstant limit 
of thick curves in $X$, then 
$h \circ \g$ is thick.
\end{cor}

\pf Let $X_n = (X_n, d/r_n)$ and suppose $(X_n, x_n) \to (T,t)$ is a tangent space.  
As in Proposition \ref{prop:tgt}, suppose $f^{k_n+n_0}:(X,d/r_n, x_n)\to (X, d, y_n)$ 
tends to $h: (T,t) \to (X,d, y)$.  
The hypothesis implies that there exist thick curves 
$\gamma_n: [0,1] \to X_n$ converging to $\gamma: [0,1] \to T$.  
The definitions of convergence of maps imply that the curve $h \circ \gamma: [0,1] \to X$ 
is the uniform limit of  the sequence of curves $\beta_n = f^{k_n+n_0}\circ \gamma_n$.  
By Proposition \ref{prop:fthick=thick}, each curve $\beta_n$ is thick.  Since nonconstant limits of thick curves are thick (Proposition \ref{prop:thick}), $\gamma$ is thick.
\qed

\section{Non-crossing curves}

In this section, we prove the following:

\begin{thm}\label{thm:noncross}
Assume that $f: S^2 \to S^2$ is topologically cxc and $X \in \GGG(f)$ is a $Q$-regular 
$2$-sphere with $Q=\confdim_{AR}\GGG(f)>2$. Then  
$f$ is topologically conjugate to a Latt\`es example.
\end{thm}

In the remainder of this section, $f$ and $X$ satisfy the hypotheses of Theorem \ref{thm:noncross}.
>From Propositions \ref{prop:positive_modulus} and \ref{prop:thick}, we obtain a family $\Gamma$ 
of thick curves of positive $Q$-modulus on $X$; by Lemma \ref{lma:cpctsup} we may assume 
it is closed and that its elements have diameters bounded above and  below.  
Since $Q>2$, Proposition \ref{prop:regular_intersection} implies that thick curves on $X$ cannot cross.

\subsection{Foliation by thick curves}

\begin{prop}\label{prop:fol1} There exist a tangent space $T$ of $X$ and 
a foliation $\FFF_T$ of $T$ by bi-infinite geodesics such that each leaf
is a limit of rescaled thick curves of $X$.
\end{prop}

>From this, we will deduce:

\begin{cor}\label{cor:fol1} The space $X$ admits a foliation $\FFF_X$ with finitely
many singularities invariant
by $f$. More precisely, there exists a finite set $F$ such that $X\setminus F$ is 
foliated  by locally thick curves, each point of $F$ is a one-prong singularity of $\FFF_X$, and 
$f^{-1}(\FFF_X\setminus F)\subset \FFF_X\setminus F$.\end{cor}

This proves Theorem \ref{thm:fol} in the iterated case.

\gap

By a one-prong singularity of a foliation of a surface, we mean that the foliation near that 
point is equivalent to the singular foliation of the plane near the origin obtained by starting with the horizontal foliation on the complex plane and identifying points $z$ and $-z$.

\gap

The proofs occupy the remainder of this subsection.

\subsubsection{Thick crosscuts}

Let $D$ be a Jordan domain in $X$. A {\em crosscut} is a curve $\g$ such that 
$\g(0),\g(1)\in\partial D$, $\g(0)\ne\g(1)$, and $\g(0,1)\subset D$. Let us say
that a crosscut $\g$ is thick if, for all $\ep>0$, $\mod_Q\G_\ep >0$, where $\G_\ep$
denotes the subset of crosscuts of $B_\infty(\g,\ep)$.

\gap

The main results of this paragraph are

\begin{prop}\label{prop:nc0} Let $\g$ be a thick crosscut of $D$.
\be
\item There are exactly two components $D_\pm$ of $D\setminus \g(0,1)$ for
which the boundary intersects $\partial D$.
\item The image $\g[0,1]$ is either contained in $\partial D_+$ or in $\partial D_-$.
\item Any other component $W$ is a Jordan domain, and there are parameters
$0<s<t<1$ such that $\partial W\subset\g([s,t])$ with $\g(s)=\g(t)$.
\eb
\end{prop}

and

\begin{prop}\label{prop:nc} Let $\g_0$ and $\g_+$ and $\g_-$ be three
thick crosscuts of $D$ with endpoints in the following cyclic order: 
$$\g_+(0)<\g_0(0)<\g_-(0)<\g_-(1)\le \g_0(1)\le \g_+(1) < \g_+(0)\,.$$ 
Then $\g_0\cap\g_-\cap\g_+\cap D=\emptyset$.
\end{prop}

We first establish some preliminary facts.
The definitions imply at once the following fact:
\begin{fact}\label{fact:nc1} Let $\g,\g'$ be two crosscuts with
endpoints in the following cyclic order  $$\g(0)<\g'(0)<\g(1)<\g'(1)\,.$$
Then $\g$ and $\g'$ cannot be thick simultaneously.\end{fact}

\begin{fact}\label{fact:nc2} Let $\g_0$ be a thick crosscut, and let $\g_1:[0,1]\to D$
be thick. Then $\g_1$ intersects at most one connected component of $D\setminus\g_0$.\end{fact}

\pf Let us proceed by contradiction and assume that $\g_1(0)$ and $\g_1(1)$ lie in different components
$U$ and $V$ of $D\setminus\g_0$. By Fact \ref{fact:nc1}, we may also assume that
$\partial U\subset \g_0$.

Let $\tau=\sup\{t>0,\ \g_1(t)\in U\}$ and $z_0=\g_1(\tau)$. By construction,
(a) $z_0\in\partial U$; (b) for any $\ep>0$, there is some $t\in (\tau-\ep,\tau)$ with 
$\g_1(t)\in U$.

Let $s\in (0,1)$ with $\g_1([s,1])\subset V$. Let $D'\subset D$ be a Jordan neighborhood of $\g_1([\tau,s])$. If it is small enough,
then $\g_1|_{[a,b]}$ is a crosscut of $D'$, where $(a,b)$ is the connected
component of $\g_1^{-1}(D')$ containing $\tau$, with $\g_1(a)\in U$ and 
$\g_1(b) \in V$. Furthermore, the connected component $\g_0'$
of $\g_0\cap D'$ which contains $z_0$ separates $\g_1(a)$ from $\g_1(b)$ otherwise we would have $U=V$.

Therefore, there are two intervals $[s_1,s_2]$ and $[t_1,t_2]$ such that
$\g_0(s_1),\g_0(t_1)\in\partial D'$, $\g_0(s_2)=\g_0(t_2)=z_0$ and
$\g_0((s_1,s_2])\cup \g_0((t_1,t_2])$ is the image of a crosscut of $D'$ which separates
$\g_1|_{[a,b]}$. Since $\g_1$ is thick, one may find another thick curve $\g_2$
arbitrarily close to $\g_1$ which avoids $z_0$. We may thus extract a crosscut
of $D'$ from $\g_2$ which either intersects $\g_0((s_1,s_2))$ or  $\g_0((t_1,t_2))$.
In both cases, we may reduce $D'$ and apply  Fact \ref{fact:nc1} to obtain a contradiction.
\qed

\gap

We now turn to the proofs of the propositions. 

\gap

\pf (Prop.\,\ref{prop:nc0})
The first point is clearly true since $\g$ is a crosscut. Note that if $\g$ is not contained in $\partial D_+$, 
then there is some parameter $t$ such that $d(\g(t),D_+)\ge \de >0$. Therefore, any thick curve at distance
at most $\de/2$ from $\g$ cannot lie in $D_+$ by Fact \ref{fact:nc2}; 
similarly for $D_-$. But any crosscut has to intersect $D_+\cup D_-$, 
so Fact \ref{fact:nc2} yields a contradiction. This proves 2.

Let us prove the last point: let $W$ be a component of $D\setminus\g(0,1)$ different from
$D_\pm$. Then $\partial W$ is contained in $\g(0,1)$ which is locally connected, so that 
Carath\'eodory's theorem implies that any conformal map $h:\D\to W$ extends continuously to
the boundary. If $\partial W$ is not a Jordan domain, then there are two rays in $\D$ which
are mapped to a Jordan curve in $\overline{W}$ which separates $\g$. Fact \ref{fact:nc2} yields another contradiction.
So $W$ is a Jordan domain. Let $s=\min \g^{-1}(\partial W)$ and $t=\max\g^{-1}(\partial W)$:
by construction, $\partial W\subset \g([s,t])$.
If $\g(s)\ne\g(t)$, then $\partial W\setminus \{\g(s),\g(t)\}$ are two arcs $a_\pm$.
Note that since $c_\pm =\g[0,s]\cup a_\pm \cup \g[t,1]$ are crosscuts contained in $\g$ and 
$a_+\cap c_-=a_-\cap c_+=\emptyset$, point 2. above cannot hold, so we obtain a contradiction
and $\g(s)=\g(t)$. \qed

\gap

\pf (Prop.\,\ref{prop:nc}) It follows from Fact \ref{fact:nc2} that we may assume that
$\g_-$ and $\g_+$ lie in different components of $D\setminus\g_0$.
Let us assume that they all meet at a point $z_0\in D$. 

One may find some $\ep>0$ such that $\g_\pm\notin B_\infty(\g_0,\ep)$.
By Fact \ref{fact:nc2}, any thick curve in $B_\infty(\g_0,\ep)$
is squeezed between $\g_-$ and $\g_+$: thus, they all go through $z_0$, which is impossible
by Proposition \ref{prop:upperbound} and the definition of thick curves. 
\qed

\subsubsection{Relation order on thick curves}
\label{subsubsecn:order}

Let $\g_0:[0,1]\to X$  be a parametrized rectifiable thick curve of $X$. 
We will define an ordering on a space $\Gamma$ of parametrized thick curves $\gamma$ that are close to 
$\g_0$ in $\PPP(X)$. This ordering will be a reflexive and transitive binary relation
(but might not be antisymmetric).

First, we choose $\g_0$ in a convenient way.  
Let $V$ be a Jordan domain which contains $\g_0$ in its interior.  
Decompose $V$ into three open Jordan subdomains $V_0$, $V_-$ and $V_+$ as shown.  That is, 
there is a homeomorphism of $V$ to $[-2,2]\times [0,1]$ such that $V_-$ is mapped to $[- 2,- 1]\times [0,1]$, 
$V_+$ is mapped to $[1,2]\times [0,1]$, $\gamma(0) \in V_-$ and $\gamma(1) \in V_+$.  
\begin{center}
\includegraphics[width=3in]{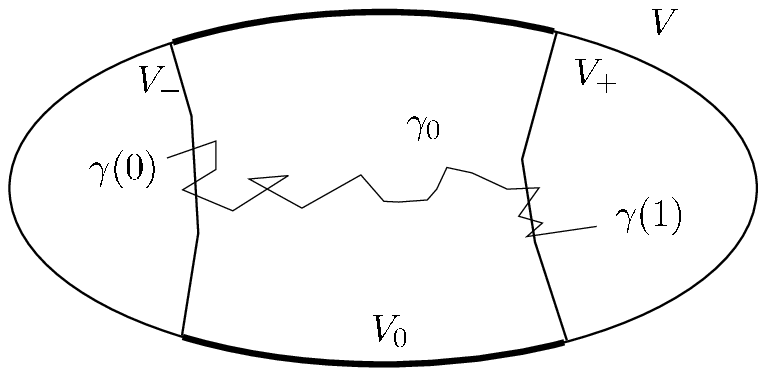}
\end{center}
Since $\gamma_0$ is rectifiable and compact, there are at most finitely many subcurves of $\gamma_0$ in $V$ 
connecting $V_-$ to $V_+$.   
Restricting to such a subcurve if necessary and reparametrizing, 
Proposition \ref{prop:thick} implies that we may assume  that there is exactly one such subcurve.  

\gap

Fix $\gamma_0, V_0, V_\pm, V$ as in the preceding paragraph.  
Next, we define $\Gamma$, a set of suitable curves close to $\gamma_0$.  There exists $r>0$ small enough so that 
for any curve $\g \in \overline{B_{\infty}(\g_0,r)}$, we have (i) $\g \subset V$, 
(ii) $\g(0) \in V_-, \g(1) \in V_+$, and (iii) there is exactly  one connected component of $\g\cap V_0$ whose 
closure 
intersects both $\partial V_+$ and $\partial V_-$; we denote this distinguished component by $\gamma_d$.  
By Proposition \ref{lma:cpctsup}, there exists 
a compact nondegenerate subfamily $\G\subset \overline{B_{\infty}(\g_0,r)}$ comprised of rectifiable thick 
curves such that $\mod_Q\G>0$. 

Finally, we define an ordering  $\leq$ on curves in $\G$.
Given $\g\in\G$, the distinguished component $\gamma_d$ cuts $V_0$ into at least two components. 
We let $U_+(\g)$ denote the component
of $V_0\setminus\g_d$  the boundary of which (when seen in the chart $V \to [-2,2]\times [0,1]$) 
contains the arc $[-1,1]\times \{1\}$, and $U_-(\g)$
the one whose boundary, when similarly viewed, contains $[-1,1]\times \{0\}$; 
these boundary arcs are indicated in bold in the figure.

The fact that thick curves do not cross implies that given $\g,\g'\in\G$, either
$U_+(\g)\subset U_+(\g')$ and $U_-(\g)\supset U_-(\g')$, or 
$U_+(\g)\supset U_+(\g')$ and $U_-(\g)\subset U_-(\g')$.
In the former case, we write $\g \geq \g'$.

\subsubsection{Definition of the tangent space}\label{subscn:def_tangent}

The tangent space will be a based at a density point $y$ of a set $A$ which we now define.
We may assume each $\gamma \in \Gamma$ is parametrized by arc length,
the parameter lying in a compact interval $I_\gamma$, and is extended so as to be completely parametrized.
We may assume $\Gamma$ contains all complete parametrizations of its curves.
For $\g\in\G$ and $t \in I_\g$
we let $$f_k(\g,t)= \sup_{|h|\le 1/k} \left(1 - \frac{d(\g(t-h),\g(t+h))}{2|h|}\right);$$
this measures how much $\gamma$ deviates from being geodesic near $t$.
Since $\g$ is rectifiable and parametrized by arclength, Proposition \ref{prop:rect} ensures that
$f_k(\g,t) \to 0$ as $k \to \infty$ for all $\g\in\G$ and almost all $t\in I_{\g}$. Set, for $x\in\supp\, \G\cap V$,
$$g_k(x)=\min \{f_k(\g,t): \gamma \in \Gamma,  \g(t)=x\}$$
and let $A$ be the set of points $x$ such that $\lim_k g_k(x)=0$.
We claim $A$ has positive measure.  Otherwise, we could define an admissible metric
$\rho$ by choosing $\rho=\infty$ on $A^*$,
and $\rho=0$ otherwise. This function is $L^Q$ integrable, and for all $\g\in\G$ and almost all $t\in I_{\g}$,
$\{f_k(\g,t)\}_k$ tends to zero, implying that $\int_\g\rho = \infty$ for all $\g\in\G$: this would
contradict  $\mod_Q\G>0$. Therefore, $A$ has positive measure.
By Egorov's theorem, we may assume that $(g_k)$ tends to zero uniformly on $A$.
We pick a point of $t$-density $y\in A$, which exists by Lemma \ref{lma:mt_density}.

The point $y$ belongs to the image of a curve $\g$.
Fix a positive sequence $(r_n)$ tending to $0$. Without loss of generality, we may
assume that $(X_n,d_n,y)=(X,d/r_n,y)$ tends to a metric plane $(T,t)$, and
that $(f^{k_n})_n$ tends to
a map $h:T\to X$ with a discrete branch set (Theorem \ref{thm:limit_surface} and Proposition \ref{prop:tgt}).

The definition of $\Gamma$, the definition of $A$, and Proposition \ref{coro:gdq}
imply that every point in $T$ belongs to the image of a geodesic
$\gamma: \R \to T$ which is a limit of a sequence $(\gamma_k)$ of rescaled thick
curves passing through points of $A$.
We denote by $\FFF_T$ the set of geodesics $\gamma$ obtained in this way. Note that $\FFF_T$ is closed with
respect to the compact-open topology.

\begin{lemma}\label{lma:T_into_2} 
With the notation from above, the geodesic $\g$ 
separates $T$ into two
simply connected regions $U_+(\g)$ and $U_-(\g)$, each of which is the interior of the respective limit 
of  $\overline{U_+(\g_k)}$ and $\overline{U_-(\g_k)}$. \end{lemma}

\pf The curve $\g$ is a simple curve, unbounded on both sides: considering the Alexandroff compactification
of $T$, $\overline{\g}$ becomes a simple closed curve so that the Jordan theorem implies that $T\setminus\g$
is the union of two disjoint open disks.

Denote by $K_\pm(\g)$ the set of points $w \in T$ arising as limits (as $k \to \infty$) 
of sequences of points $w_k \in U_\pm(\g_k)$, respectively. 
Let $w\in K_+(\g)\cap K_-(\g)$, and $(w_k)$, $(w_k')$ be sequences as above which tend to $w$;
by the bounded turning (BT) property, we may find
continua $C_k\subset X_{n_k}$ joining $w_k$ and $w_k'$ with
$\diam\,C_k\lesssim |w_k-w_k'|$. By the Jordan curve theorem,
it follows that $C_k\cap\g_k\ne\emptyset$. As $k$ tends to infinity,
we obtain that $w\in\g$. Hence $K_+(\g)\cap K_-(\g)\subset\g$. 

Let us prove that $K_+(\g)\cup K_-(\g)=T$. Pick $w\in T$, and $(w_k)\in\prod X_{n_k}$ which tends
to $w$. Either, there are infinitely many $k$ such that $w_k\in \overline{U_+(\g_k)\cup U_-(\g_k)}$, in which
case $w\in K_+(\g)\cup K_-(\g)$. Or, $w_k$ is contained in a bounded component $W$ of $X_{n_k}\setminus\g_k$
for all $k$ large enough. 

In the latter case, we apply Proposition \ref{prop:nc0} to obtain $s<t$
such that $\g_k(s)=\g_k(t)$, $\g_k[s,t]$ separates $w_k$ from $U_+(\g_k)\cup U_-(\g_k)$
and $\g_k(s)\in\partial (U_+(\g_k)\cup U_-(\g_k))$.

Since $\g_k(s)=\g_k(t)$, we have $$|s-t|=|\g(s)-\g(t)|\le 2\|\g_k-\g\|_\infty$$
so that $\diam\, \g_k([s,t])\le 2\|\g_k-\g\|_\infty$. Property (ALC2) now implies that
$$\diam W\le 2L\diam \g_k([s,t])\le 4 L \|\g_k-\g\|_\infty\,.$$

It follows that we may find $w_k'\in\overline{U_+(\g_k)\cup U_-(\g_k)}$ with 
$$|w_k-w_k'|\lesssim \|\g-\g_k\|_\infty\,.$$
This proves that $w\in K_+(\g)\cup K_-(\g)$ so $K_+(\g)\cup K_-(\g)=T$. The proof follows.\qed

\begin{remark} If $\g,\g'\in\FFF_T$, then by Corollary \ref{cor:locallythick} they cannot cross 
at a regular point of branched covering $h: T \to X$.\end{remark}

\subsubsection{Simultaneous zooms of thick curves}

In this section and the next, we examine the structure of the set of geodesics $\FFF_T$ in the metric plane $T$ constructed 
in the previous subsection.  
Recall that the geodesics in $\FFF_T$ were obtained as limits of rescaled curves in the family $\Gamma$ defined in \S\,\ref{subsubsecn:order}.  

Let $\g,\g'\in\FFF_T$. We shall write $\g\le \g'$ if there exist a subsequence $(n_k)_k$ and sequences
of curves $(\g_k)$ and $(\g_k')$  such that $\g_k,\g_k'\subset X_{n_k}$, $\g_k\le \g_k'$ and if
both sequences converge to $\g$ and $\g'$ respectively. 

We now fix $\g,\g'\in\FFF_T$ and assume that there exist  a subsequence $(n_k)_k$ and sequences
of curves $(\g_k)$ and $(\g_k')$  with $\g_k,\g_k'\subset X_{n_k}$ and 
such that both sequences converge to $\g$ and $\g'$ respectively. For each $k$, 
either $\g_k \leq \g_k'$ or $\g_k' \leq \g_k$. 
So, extracting a subsequence if necessary, we may  assume
that $\g\le \g'$. 
We may also assume that $\overline{U_\pm(\g_k)}$ and  $\overline{U_\pm(\g_k')}$ 
converge to closed  and connected domains $K_\pm(\g)$ and $K_\pm(\g')$, cf. Lemma \ref{lma:T_into_2}. 

\begin{lemma}\label{lma:order1} If $\g\ne\g'$, then, for any $z\in U_+(\g)\cap U_-(\g')$,
there exists a curve $\g_z\in\FFF_T$ going through $z$ such that $\g\le\g_z\le
\g'$.\end{lemma}

\pf Suppose $z_n \to z$, $z_n \in X_n$. We may assume that
$z_n\in U_+(\g_n)\cap U_-(\g_n')$. By Lemma \ref{lma:density_tangent}, there is a sequence
of curves $(\g_{z,n})_n$ such that $d_n(z_n,\g_{z,n})$ tends to $0$.
It follows that $\g_n\le \g_{z,n}\le\g_n'$. Extracting a limit, the
lemma follows.\qed

\begin{lemma}\label{lma:inter1}
If $\g$ and $\g'$ intersect, then $\g=\g'$.\end{lemma}

\pf We assume that $\g$ and $\g'$ are distinct but that $\g(0)=\g'(0)=z_0$
and $\g((-a,0))\cap\g'=\g'((-a,0))\cap\g=\emptyset$.

Let $z\in U_+(\g)\cap U_-(\g')$.  By Lemma \ref{lma:order1}, 
there exists a curve $\g_z\in\FFF_T$ going through $z$ such that $\g\le\g_z\le
\g'$ For any $n$, there exists a continuum
$C_n$ joining $\g_n(0)$ to $\g_n'(0)$ with $\diam\,C_n\lesssim |\g_n(0)-\g_n'(0)|$. 
By the Jordan curve theorem, $C_n\cap\g_{z,n}\ne\emptyset$, so
that $z_0\in\g_z$. 

We note that $z_0$ belongs to three limits of thick curves.
If $h$ is a local homeomorphism at $z_0$, then we obtain a contradiction by Proposition \ref{prop:nc}.

We assume now that $h$ is not a local homeomorphism at $z_0$.  Set $w_0 = h(z_0)$.  
By Corollary \ref{cor:tangent_map_surface}, 
we may assume that in local coordinates near $z_0$ and $w_0$, the map $h$ is given by 
$h:\D\to\D$, $h(z)=z^k$, for some $k\ge 2$; in particular, this restriction is proper.  For each element $c\in\FFF_T$ which contains the origin,
we may restrict $c$ and reparametrize it so that $c:[-1,1]\to\overline{\D}$ is a crosscut
with $c(0)=0$.

We now establish some facts about crosscuts $c$ obtained as the restrictions
of elements $C$ of $\FFF_T$ with $\g\le C\le \g'$.
\be
\item Since $h$ is proper on $\D$, $h(c)$ joins $0$ to the boundary of $\D$.
\item Corollary \ref{cor:locallythick} implies $h(c)$ is also thick, so $h(c)\cap h(\g)=\{0\}$ since both images are thick curves.
\item The map $h$ is locally injective on $\D\setminus \{0\}$ and $h^{-1}(0)=0$, so, if $(h\circ c)|_{[-1,0]}$ was not injective,  
then
there would be different times $s<t$ in $(-1,0)$ such that $h(c(s))=h(c(t))$.  But then  the thick
curves $c|_{[-1, (2s+t)/3]}$ and $c|_{[(s+2t)/3,0]}$ intersect at a regular point of $h$, which is impossible.
Similarly, $(h\circ c)|_{[0,1]}$ is also injective. It also follows that if $h\circ c$ is not 
globally injective, then $c([-1,0])=c([0,1])$.
\item Furthermore, there may be at most $(k-1)$ ordered curves
which are mapped 2-to-1 under $h|_{\D}$, so that we may as well consider that
$h(c\cap \D)$ is an arc for every restriction of a curve $c\in\FFF_T$ which defines
a crosscut and which
goes through the origin such that $\g\le c\le \g'$. 
\eb

Given the curve $\g$ in the lemma, we let $\Omega$ be the connected component 
of $\D\setminus h^{-1}(h(\g([-1,0])))$
which contains $\g$ in its boundary and intersects $U_+(\g)$. 
Note that $h|_{\Omega}:\Omega\to \D\setminus h(\g((-1,0]))$ is a 
homeomorphism.
For any other curve $c\in\FFF_T$, we know that if $h(c)$ is not a subset of $h(\g)$, then $c$ cannot
intersect $h^{-1}(h(\g([-1,0))))$; therefore, we may define $\widehat{c}\subset \Omega\cup\g([-1,0])$ as the lift of $h(c)$
to $\Omega\cup\g([-1,0])$.

It follows from above that $\widehat{\g'}\cap \Omega$ is contained in a component of $\Omega\setminus\widehat{\g}$, 
for otherwise $h(\g')$ and $h(\g)$ would cross. 
Let $V'$
be the connected component of $\Omega\setminus\widehat{\g'}$ with boundary $\widehat{\g'}$.
We first consider the case that $V'$ is not a subset of $U_+(\g')$. Then, 
for any $z \in V'\setminus U_+(\g')$, Lemma \ref{lma:order1} provides us with a curve
$\g_z\in\FFF_T$ such that $\g\le \g_z\le \g'$. In particular $0\in\g_z$. 
Since $\g_z\cap V'\ne\emptyset$, it follows that $\widehat{\g_z}\subset V'\cup\{0\}$ so that
$h(\g)$, $h(\g')$ and $h(\g_z)$ satisfy the assumptions of Proposition \ref{prop:nc}, 
yielding a contradiction.
Therefore, $V'\subset U_+(\g')$. We consider a new curve $\g''$ given by Lemma \ref{lma:order1}
such that $\g\le \g''\le \g'$ and $\g''\cap \D$ is not contained in $U_+(\g')$: let $V''$ be
the connected component of $\Omega\setminus\widehat{\g''}$ with boundary $\widehat{\g''}$.
By construction, $V''$ is not a subset of $U_+(\g')$. As above, we consider a curve
$\g_z\in\FFF_T$ such that $\g\le \g_z\le \g'$ and $\widehat{\g_z}\subset V''\cup\{0\}$: we obtain
another contradiction from $h(\g)$, $h(\g_z)$ and $h(\g'')$.
\qed

\subsubsection{Different zooms}

Let us assume that $\g,\g'\in \FFF_T$ are two distinct curves.

This implies that there are
two increasing functions $\varphi,\psi:\N\to\N$, and thick curves 
$\g_{\varphi(n)}\subset X_{\varphi(n)}$ and $\g_{\psi(n)}'\subset X_{\psi(n)}$
which tend to $\g$ and $\g'$ respectively.

\begin{lemma} The curves $\g$ and $\g'$ are disjoint.\end{lemma}

\pf 
By the previous section, we may assume that the sequences 
$\varphi$ and $\psi$ have no integers in common.  
We will show that if $\g$ and $\g'$ are not disjoint, then we can build a curve $\g''$ which will cross $\g'$ at a regular
point of $h$, contradicting 
Corollary \ref{cor:locallythick}.

We may assume that $\g(0)=\g'(0)$. Since $h$ has a discrete
branch set (Corollary \ref{cor:tangent_map_surface}), we may assume that $\g'|_{(0,a)}$ has no branch points
of $h$, and that it is disjoint from $\g$. Let $z_1=\g'(a/2)$,
and apply Lemma \ref{lma:order1} to obtain a curve $\g_1$ going through $z_1$
obtained by blowing up thick curves from $X_{\varphi(n)}$;
we may assume that $\g\le \g_1$. 
If $\g_1$ does not cross $\g'$, then we construct a similar
curve $\g_2\,(\ne\g_1)$ intersecting $\g'$ such that $\g\le \g_2\le \g_1$.
It follows that $\g_2$ has to cross $\g'$, since it is disjoint from $\g$ and $\g_1$
by Lemma \ref{lma:inter1}.\qed

\subsubsection{Foliations}

We have proved that any point in $T$ belonged to a unique
element of $\FFF_T$, and that no two such curves could intersect.

It remains to prove that $\FFF_T$ is a genuine foliation. Following Whitney \cite[Part II]{whitney:foliationcurves},
one has to check that this is a {\it regular family of curves}:
given any point $p$ and a direction on the curve through $p$, there is an arc $pq$ 
in this direction with the property that for every $\epsilon>0$ there is a $\delta>0$ such that, 
for any point $p'$ within a distance $\delta$ of $p$, there is an arc $p'q'$ of the curve through $p'$
which lies within an $\epsilon$-neighborhood of $pq$ and on which $q'$ lies within an $\epsilon$-neighborhood of $q$;
moreover, if $r'$ and $s'$ are two points on $p'q'$ within a distance $\delta$ of each other, 
then the diameter of the arc $r's'$ is less than $\epsilon$ (see also \cite[Lme\,4.7]{kolev:sousgroupe_compact} for
an explicit construction of a transversal arc). 

The second condition is automatically satisfied since curves are geodesics of $T$.
For the first, assume that it does not hold: then, given an arc from $\FFF_T$ joining two points
$p$ and $q$, there is some $\epsilon_0>0$ with the property that there is a sequence of points
$(p_n)$ tending to $p$ such that no subarc $(p_n q_n)$ either leaves the $\epsilon_0$-neighborhood
of $(pq)$ or $q_n$ is $\epsilon_0$-away from $q$. If we consider parametrizations of the curves
$\g_n$ going through $p_n$ with $\g_n(0)=p_n$, then, Ascoli's theorem implies that 
the sequence $(\g_n)$ tends uniformly on compact sets to a geodesic $\g$  going through $p$;
since $\FFF_T$ is closed, it follows that $\g\in\FFF_T$,
contradicting the definition of $(p_n)$. Therefore $\FFF_T$ is a regular family of curves,
hence a genuine foliation. 

This ends the proof of Proposition \ref{prop:fol1}.

\gap

Before proving Corollary \ref{cor:fol1}, we note the following:

\begin{lemma}
\label{lemma:fol_by_thick}
Suppose $\FFF_1, \FFF_2$ are two nonsingular foliations of an open subset $U$ of $X$, 
each of whose leaves are locally thick curves.  Then $\FFF_1 = \FFF_2$.
\end{lemma}

\pf Let $x \in U$.  The leaves $L_1, L_2$ in $\FFF_1, \FFF_2$ containing $x$ cannot cross 
at $x$ since they are both thick. 
If they are not the same, there is a leaf $L_1'$ near $L_1$ which crosses $L_2$, 
which is impossible.  
\qed

\pf (Cor.\,\ref{cor:fol1}) We push the foliation $\FFF_T$ down to a singular foliation 
$\FFF_X$ of $X$ using $h$: away from the images of branch points of $h$, 
we define the leaves of $\FFF_X$ 
to be the images of the leaves of $\FFF_T$ under $h${
by Corollary  \ref{cor:locallythick}, they are  locally thick curves, so 
Lemma \ref{lemma:fol_by_thick} implies $\FFF_X$ is well-defined.

The construction of the tangent $T$ implies that there exists a compact set $K\subset T$ 
such that $h(int(K))=X$.  Since $h$ is a branched covering, its branch locus $B_h$ is discrete, 
hence its intersection with $K$ is finite. For any regular point $z$ of $h$ in $int(K)$,
$h_*\FFF_T=\FFF_X$ is regular in the neighborhood of its image $h(z)$.
It follows that the set $F$ of singular points of $\FFF_X$ is contained in the image 
$h(int(K) \intersect B_h)$ and is therefore finite.

We now analyze the structure of the singularities of $\FFF_X$.  
Let $x_0$ be a branch point of $h$; we use the fact that $T$ and $X$ are surfaces and $h$ is an open map: 
by Corollary \ref{cor:tangent_map_surface}, there are neighborhoods $U$ and $V$ 
of $x_0$ and $h(x_0)$ respectively such that $h:(U,x_0)\to (V,h(x_0))$ 
is equivalent to $z\mapsto z^{k}$  in the unit disk of the complex plane $\C$ 
for some positive integer $k\ge 2$. 
Let $\g$ be the connected component in $U$ of the leaf of $\FFF_T$ which contains $x_0$;
since $h|_U$ is proper onto $V$, there is some $j\in\{1,2\}$ such that its image $h(\g)$ cuts $V$ into $j$ components,
so that $h^{-1}(h(\g))$ cuts $U$ into $j\cdot k$ components. Since thick curves on $X$ cannot cross, it follows
that $h^{-1}(h(\g))\cap U=\g$ and $j\cdot k=2$: hence $k=2$ and $j=1$. Therefore, the branched covering 
$h:T\to X$ can have only simple critical points 
and a curve which goes through a critical point of $h$
must map in a locally 2-to-1 fashion onto its image: all singularities of $\FFF_X$ are prongs.

We now prove the invariance of the foliation $\FFF_X$ under $f$: if $x\in X$ is not a branched point
of $f$, then $x\in F$ if and only if $f(x)\in F$. Furthermore, if $f(x)$ is a regular point of $\FFF_X$,
then, pulling back $\FFF_X$ under $f^{-1}$ to a neighborhood of $x$, we see that $x$ cannot be a branched
point for $\FFF_X$ would have a singularity which would not be a prong: the invariance of $\FFF_X$ is established.

This ends the proof of Corollary \ref{cor:fol1}.
\qed

\subsection{Parabolic orbifold structure}\label{subscn:parorbi}
In this subsection, we prove that $f$ is conjugate to a Latt\`es example. We will rely on
the ``easy part'' of Douady and Hubbard's classification of finite branched coverings of the sphere \cite{DH1}.

\gap

\pf (Theorem \ref{thm:noncross})
By Corollary \ref{cor:fol1}, we know that $f^{-1}(\FFF_X\setminus F)\subset \FFF_X\setminus F$.
It follows that  $f^{-1}(F)\subset (B_f\cup F)$, so that $P_f\subset F$
and $f$ is postcritically finite. 
Note that if $f$ is locally injective at $x$, and if $y=f(x) \in F$, then $x \in F$.  
This implies that in fact $P_f = F$, that $f^{-1}(P_f) = P_f \union B_f$.  
>From \cite[Lma\,3.2]{DH1} it follows that $\#P_f \leq 4$ and more precisely that in fact $\#P_f=4$.
Since all the singularities of $\FFF_X$ are prongs, the branch points of $f$ and of $h$ are simple, 
so the orbifold associated to $f$ is the $(2,2,2,2)$-orbifold.
We may now apply \cite[Prop.\,9.3]{DH1} to deduce that  
$f$  lifts as a covering map $g$
of a torus $\tilde{X}$. 
Since $f$ is topologically cxc, $g$ is positively expansive.  Pick a basis for $H_1(\tilde{X}, \Z)$ and let $A$ be the matrix of the induced map 
$H_1(g): H_1(\tilde{X}, \Z) \to H_1(\tilde{X}, \Z)$.  Since $g$ is positively expansive, $g$ and the map on the torus $\R^2/\Z^2$ induced 
by $A$ are topologically conjugate.  Hence $f$ is topologically conjugate to the Latt\`es example $f_A$ induced by $A$. 
\qed

\subsection{Tangents as universal orbifold covering}

In this subsection, we establish the claim in  Remark \ref{rmk:orbifold-cover}.

We begin with an easy consequence of the analysis in the preceding subsection: 

\begin{cor}
\label{cor:univcover}
Under the hypothesis of Theorem \ref{thm:noncross}, the tangent map $h: T \to X$ constructed 
in \S\,\ref{subscn:def_tangent}  is the universal orbifold covering map of the orbifold associated to $f$.
\end{cor}

\pf Since thick curves cannot cross, every branch point of $h$ is simple and maps to an element of $F$ and, conversely, any preimage of a point in $F$ under $h$ is a simple branch point of $h$.  
\qed

\gap

The construction of the above tangent space $T$ is very indirect.  We continue 
with a proposition which gives a direct construction in which the basepoint is a periodic point.

\begin{prop}\label{prop:tgtrepelling} Let $f:S^2\to S^2$ be a topological 
cxc map, $p$ be a periodic point of $f$ of period $k$, and suppose there is a neighborhood $W$ of $p$ such that 
$f^k: W \to f^k(W)$ is a homeomorphism with $\cl{W} \subset f^k(W)$. Let $X\in\GGG(f)$.

Then there exists a tangent space $(T,t)$ to $X$ at $p$ locally homeomorphic to $W$, an 
expanding homeomorphism $\psi:(T,t)\to (T,t)$ fixing $t$ whose iterates are uniformly quasisymmetric, 
and a quasiregular map $h:(T,t)\to (X,p)$ such that 
$h\circ\psi= f^k\circ h$, $h(t)=p$.  

In the case of the sphere equipped with a visual metric, i.e. $X=(S^2, d_v)$, the map $h$ is an isometry near $t$, 
and the iterates of $\psi$ are similarities with constant expansion factor.  \end{prop}

Here, by quasiregular, we mean an open map which is locally uniformly quasisymmetric
away from its discrete branch set.

\gap

\pf By Lemma \ref{lma:tgt_equiv}, it is enough to prove the proposition with $d=d_v$.

Let $g:(W,p)\to (f^{-k}(W),p)$ be the local inverse of $f^k$ near $p$. We may
assume that $|g(x)-g(y)|=\lambda |x-y|$ with $\lambda =\te^k$,
for all $x,y\in W$ by \cite[Prop.\,3.2.3]{kmp:ph:cxci}; here, $\theta$ is the constant given by Theorem \ref{thm:can_gauge}.

Define scaling functions $\sigma_n:(W,d)\to W_n=(W, d/\lambda^n)$
and let $g_n:W_n\to W_{n+1}$ be defined as $g_n=\sigma_{n+1}\circ g\circ \sigma_n^{-1}$, which is an isometric embedding.

We may then consider the following inductive limit:
$$\WWW=\lim_{\longrightarrow}(W_n,g_n)\,.$$

Since the $g_n$ are isometries, the set $\WWW$ is naturally a metric space.
For all $n$, $W_n$ embeds canonically in $\WWW$. 
Define $\psi:\WWW\to\WWW$ by $\psi_n(x)=\sigma_{n+1}\circ\sigma_n^{-1}(x)$, for $x\in W_n$. We check that
\begin{eqnarray*}
\psi_{n+1}\circ g_n & = & \sigma_{n+2}\circ \sigma_{n+1}^{-1}\circ \sigma_{n+1}\circ g\circ \sigma_n^{-1}\\
& = & \sigma_{n+2}\circ g\circ \sigma_{n+1}^{-1}\circ \sigma_{n+1}\circ \sigma_n^{-1}\\
& = &  g_{n+1}\circ\psi_n\,.\end{eqnarray*}
It follows that 
$$|\psi(x)-\psi(y)|=\frac{1}{\lambda}|x-y|\,.$$

If $x\in W_n$,
let $h_n(x)=f^{kn}\circ \sigma_n^{-1}(x)$. This defines a $1$-Lipschitz
map $h:\WWW\to X$ since
\begin{eqnarray*}
h_{n+1}\circ g_n & = & f^{kn}\circ f^k\circ \sigma_{n+1}^{-1}\circ \sigma_{n+1}\circ g\circ \sigma_n^{-1}\\
& = & f^{kn}\circ (f^k\circ g) \circ  \sigma_n^{-1}= h_n\,.\end{eqnarray*}

One has 
\begin{eqnarray*}
h_{n+1}\circ \psi_{n}  & = & f^{kn+k}\circ \sigma_{n+1}^{-1}\circ \sigma_{n+1}\circ \sigma_n^{-1}\\
& = & f^k\circ ( f^{kn}\circ  \sigma_n^{-1})=  f^k \circ h_n(x)\,.\end{eqnarray*}
Therefore, $h\circ \psi= f^k \circ h$. 

Since the maps $$\sigma_n:\left(W,\frac{d}{\lambda^n},p\right)\to (W_n,p)\subset \WWW$$ are isometries, it follows that $\WWW$ is the Gromov-Hausdorff
limit of $(X,d_v/\lambda^n,p)$.  
Note that if $W=W_0$ is chosen sufficiently small, then $h|_{W_1}$ is an isometry, by construction.  \qed

\begin{remark} If $d=d_v$,  one of the visual metrics given by Theorem \ref{thm:can_gauge}, then the assumption [Deg] in Proposition \ref{prop:tgtrepelling} can be omitted,
provided that the branch set is disjoint from the cycle containing $p$.
\end{remark}

We conclude this section with the proof of the claim in Remark \ref{rmk:orbifold-cover}.  
The Lefschetz formula shows that $f$ has a fixed-point, $p$.
Let $(T,t)$ be the tangent space at $p$, $h: (T,t) \to (X, p)$ be the quasiregular map, and 
$\psi: (T,t) \to (T,t)$  be the associated expanding homothety, each given by 
Proposition \ref{prop:tgtrepelling}. 
Let $\FFF_X$ be the foliation of $X$ by thick curves given by Corollary \ref{cor:fol1}. 

Restricting the neighborhood $W$ if necessary, we may assume
that $\FFF_X\cap W$ has at most one singularity ---which would be at $p$--- 
and that there is a neighborhood $N \subset T$ of $t$ 
such that $h:N\to W$ is quasisymmetric. The foliation 
$\FFF_N=h^{-1}(\FFF_X)\cap N$ has at most one singularity ---at $t$--- and all the leaves of $\FFF_N$ are 
thick by Proposition \ref{prop:tyson}.
Applying iterates of $\psi$, we obtain a $\psi$-invariant foliation $\FFF_T$ of $T$ by locally thick curves 
with
at most one singularity such that $h:(T,\FFF_T)\to (X,\FFF_X)$, since $\cup_n \psi^n(N)=T$.

If $\FFF_T$ has no singularity at $t$, then the same analysis as in Corollary \ref{cor:univcover} 
shows that $h:T\to X$ is again a universal orbifold cover.  
In the notation of Remark \ref{rmk:orbifold-cover}, we set $\tilde{X}=T$, $\pi=h$, and $\psi=\psi$.  
If $t$ is a singularity of $\FFF_T$, let $q: (\tilde{T}, \tilde{t}) \to (T,t)$ be a double-cover of $T$ 
ramified above $t$ (which exists, since $T$ is a plane).  
The metric on $T$ lifts to $\tilde{T}$ so that $q$ becomes a local isometry away from $\tilde{t}$.  
Then $h \circ q: \tilde{T} \to X$ is now the universal orbifold cover, and $\psi$ lifts to 
$\tilde{\psi}: \tilde{T} \to \tilde{T}$.  We set $\tilde{X}=\tilde{T}$ , 
$\pi = h \circ q$, and $\psi=\tilde{\psi}$.  
\qed
\section{Latt\`es or rational}

This section is devoted to the proof of Theorems \ref{thm:rl} and \ref{thm:fol} in the case of iterated maps. 
We let $f:S^2\to S^2$
be topologically cxc and we assume that there is an Ahlfors-regular metric space $X \in \GGG(f)$ of 
dimension $Q=\confdim_{AR}(f)$.
>From Propositions \ref{prop:positive_modulus} and \ref{prop:thick}, we obtain a family $\Gamma$ 
of thick curves of positive $Q$-modulus on $X$.

\gap

If there are two thick curves which cross, then Proposition \ref{prop:mainthmrat} implies
that $f$ is conjugate to a rational map. By \cite[Cor.\,4.4.2]{kmp:ph:cxci}, this rational map
is semihyperbolic.

\gap

If thick curves do not cross, then Corollary \ref{cor:fol1} implies the existence
of an invariant foliation ---completing the proof of Theorem \ref{thm:fol}--- and
Theorem \ref{thm:noncross} implies that
$f$ is conjugate to a Latt\`es map $f_A$.

Since $f_A$ is expanding,  the eigenvalues of $A$ have modulus larger than one.   
Since the conformal gauge is an invariant of topological conjugacy \cite[\S\,2.8]{kmp:ph:cxci},  
the conformal dimensions of $f$ and of $f_A$ are the same.  
We cannot have $Q=2$, otherwise there would be crossing thick curves.
Therefore, $Q>2$ and Theorem \ref{thm:lattesclass} implies that $A$ has two 
distinct real eigenvalues.  \qed

\appendix

\section{Further applications of the comparison formula}
Here, we sketch some applications of Proposition \ref{prop:comb_dim_estimate}
and of its corollaries. For brevity, we refer the reader to the cited references for the relevant definitions
and background.

\subsection{The Loewner property and its combinatorial version}

Following J.\,Heinonen and P.\,Koskela \cite{heinonen:koskela:loewner}, 
we say that an arcwise connected
metric space is $Q$-Loewner, $Q>1$, if there is a non-increasing function
$\psi:\R_+\to\R_+$ such that
$$\mod_Q(E,F)\ge \psi(\Delta(E,F))\,,$$
for any pair of disjoint continua $E$ and $F$ in $X$ and with
$$\Delta(E,F)=\frac{\dist(E,F)}{\min\{\diam E,\diam F\}}\,.$$ 

\gap

In \cite{kleiner:icm2006}, B.\,Kleiner suggests a combinatorial version 
of a $Q$-Ahlfors regular and $Q$-Loewner space. Let $X$ be a proper
metric space. Let $\RRR_n$ 
denote a maximal $2^{-n}$-separated set of $X$, 
then $\SSS_n=\{B(x,1/2^n), x\in\RRR_n\}$
defines a uniform sequence of quasipackings.

Say $X$ satisfies the {\em combinatorial $Q$-Loewner property} if
there are  non-increasing positive functions $\psi$ and $\varphi$ such that

$$(CLP) \quad\quad \psi(\Delta(E,F))\le \mod_Q(E,F,\SSS_n)\le \varphi(\Delta(E,F))\,,$$
for any pair of disjoint continua $E$ and $F$ in $X$ and for any $n$
large enough (with respect to $(E,F)$).

\gap

The following notion also appears in \cite{kleiner:icm2006}:
a  compact metric space  $X$ is {\em selfsimilar} if there is a constant
$L_0\ge 1$ such that for any ball $B(x,r)\subset X$ with $r\in (0,\diam\,X]$,
there is an open set $U\subset X$ which is $L_0$-bi-Lipschitz to the rescaled ball
$(B(x,r),(1/r)d)$. 

In  \cite{bourdon:kleiner:coxeter}, M.\,Bourdon and 
B.\,Kleiner prove that, for a selfsimilar space $X$, for $\de>0$ small enough, 
and for
any $p\ge 1$, 
\be
\item if $\{\mod_p(\G_\de,\SSS_n)\}_n$ is bounded, then there exists a decreasing function $\phi$ such that 
the upper bound of (CLP) holds;
\item  the sequence $\{\mod_p(\G_\de,\SSS_n)\}_n$ is essentially submultiplicative;
\item there is a critical dimension $Q_M>0$ such that $\{\mod_p(\G_\de,\SSS_n)\}_n$
tends to $0$ for $p>Q_M$, is unbounded for $1\le p<Q_M$ and admits a positive
lower bound for $p\in [1,Q_M]$.
\eb

The question arises whether the dependence on $p$ behaves like the
Hausdorff dimension. If we assume that there exists $Q>1$ such that
$$\mod_Q(\G_\de,\SSS_n)\asymp 1\,,$$
e.g. if $X$ satisfies the combinatorial Loewner property,
then $p\mapsto \{\mod_p(\G_\de,\SSS_n)\}_n$ has such a behavior according to 
Corollary \ref{cor:comb_dim}:
the moduli tend to infinity for $p<Q$ and to zero for $p>Q$ (and $Q=Q_M$).

\gap

Let $X$ be the boundary of a hyperbolic Coxeter group endowed with a visual metric.
For positive $\de$ and $r$ small enough, M.\,Bourdon and 
B.\,Kleiner define the family $\FFF^g\subset\G_\de$ of so-called generic curves
i.e., this is the subfamily of $\G_\de$ such that 
none of these curves is contained in the $r$-neighborhood of a sub-Coxeter, also called a parabolic,
subgroup. 
They introduce the critical exponent 
$$Q_m=\sup\{p\ge 1,\ \lim\mod_p(\FFF^g,\SSS_n)=+\infty\}\,.$$
If $X$  satisfies the combinatorial Loewner property,
then $Q_m=Q_M$: this follows from choosing two disjoint continua $\{E,F\}$
such that the family of curves joining them is contained in $\FFF^g$.
This provides an answer to a question which was raised 
by M.\,Bourdon and B.\,Kleiner; see \cite[Rmk\,(1), \S\,3]{bourdon:kleiner:coxeter}.

\subsection{Combinatorial moduli on surfaces}

The proof of Theorem \ref{thm:char_of_rm} relies heavily on
\cite[Thm.\,1.1]{bonk:kleiner:qsparam}:
\begin{thm}[M.\,Bonk \& B.\,Kleiner]
\label{thm:qsparam}
Suppose $X$ is a metric space which is homeomorphic to $\IS^2$, linearly locally connected, and Ahlfors $2$-regular. 
Then $X$ is quasisymmetrically equivalent  to $\IS^2$.\end{thm}

Note that if $X$ is a surface which satisfies the combinatorial
Loewner property in some dimension $Q\ge 2$, then it is 
fairly easy to construct two families of curves $\G_1$ and $\G_2$ of positive
modulus such that every curve in $\G_1$ intersects any curve in $\G_2$:
it follows at once that $Q=2$.

In particular, if $X$ is $Q$-Loewner and $Q$-regular then $Q=2$.
Hence, if $X$ is furthermore homeomorphic to $\IS^2$, then it is also
quasisymmetric to $\IS^2$: this provides an alternative proof of  
\cite[Thm.\,1.2]{bonk:kleiner:qsparam}.

\subsection{Conformal dimension and Gromov hyperbolic groups with $2$-sphere boundary}

We now sketch how  Theorem \ref{thm:hgcannon}
may be proved using Theorem \ref{thm:fol}.

The arguments in \cite{bonk:kleiner:conf_dim} used to prove Theorem \ref{thm:hgcannon} proceed by first showing that there 
exists an Ahlfors-regular Loewner metric in the gauge of the group, and then 
applying the uniformisation theorem for Loewner spheres \cite{bonk:kleiner:qsparam}.

Here, we sketch how to bypass the Loewner property and prove instead that
the boundary admits a $2$-regular metric. 

\gap

\pf (Thm \ref{thm:hgcannon}) The first step is to prove that there is a family of
positive modulus on the boundary of $G$: the argument is the same
as above. We use \cite[Lma\,5.3]{bonk:kleiner:rigidity}
which says that any weak tangent is quasi-M\"obius equivalent, say via a map $h$, to
$\partial G$, punctured at some point, cf. \cite[Cor.\,1.6]{bonk:kleiner:conf_dim}.

Then it is enough to prove that there are two thick curves which cross.
For this, we may assume that this is not the case: the argument for the proof
of  Proposition \ref{prop:fol1} works the same: the same
blowup strategy defines a weak tangent space foliated by
locally thick geodesic curves. We push this foliation forward via $h$ to $\partial G$.
Since the group acts by quasi-M\"obius homeomorphisms, it preserves
the set of thick curves, hence the foliation; this in turn implies that 
$G$ has a global fixed point (the puncture): contradiction.
Thus, there are two thick curves which cross, which implies that $Q=2$.
An application of Theorem \ref{thm:char_of_rm} concludes the proof.
\qed

\section{Applications to general topologically cxc maps}

In \cite{kmp:ph:cxci}, topologically cxc maps are defined in a much general setting than
on spheres. We recall briefly the definition and record the results which hold in this more general context.

\gap

Suppose $X, Y$ are locally compact Hausdorff spaces, and let $f: X \to Y$ be a finite-to-one continuous map.  
The {\em degree}
of $f$ is  
\[ \deg(f) = \sup\{\# f^{-1}(y): y \in Y\}.\] For $x \in X$, the {\em local degree} 
of $f$ at $x$ is
\[ \deg(f; x) = \inf_U \sup\{\#f^{-1}(\{z\}) \intersect U : z \in f(U)\}\] where $U$
ranges over all neighborhoods of $x$.   

\gap
The map $f:X\to Y$ is a {\em finite branched covering}
(abbreviated fbc) provided $\deg(f)<\infty$ and 
\bi 
\item[(i)] $$
\sum_{x \in f^{-1}(y)} \deg(f;x) = \deg f$$ holds for each $y \in Y$;
\item[(ii)] for every $x_0 \in X$, there are compact neighborhoods
$U$ and $V$ of $x_0$ and $f(x_0)$ respectively such that
\[ \sum_{x \in U, f(x)=y} \deg(f; x) = \deg(f; x_0)\]
for all $y\in V$.\ib  

\gap

Let $\XX_0, \XX_1$ be Hausdorff locally compact, locally connected topological spaces, 
each with finitely many connected components.  
We further assume that  $\XX_1$ is an open subset of $\XX_0$ and that $\cl{\XX_1}$ is compact in $\XX_0$.    
The {\em repellor} of $f: \XX_1 \to \XX_0$ is
\[ X = \{ x \in \XX_1 | f^n(x) \in \XX_1 \;\forall n > 0\} = \bigcap_n \cl{f^{-n}\XX_0}.\]

Let $\UUU_0$ be a finite cover of $X$ by open, connected subsets of $\XX_1$ whose intersection with $X$ is nonempty.  
Let $\UUU_n$, $n\ge 1$, be the connected components of $f^{-n}(U)$,  where $U$ ranges over $\UUU_0$.

We say $f: (\XX_1,X) \to (\XX_0, X)$ is {\em topologically coarse expanding conformal with repellor $X$}
provided 
\be
\item  the restriction $f|X: X \to X$ is also an fbc of degree equal to $d$;
\item there exists  a finite covering $\UUU_0$ as above, such that the axioms [Exp], [Irred] and [Deg] hold.
\eb

Without any changes, Theorem \ref{thm:can_gauge}, Proposition \ref{prop:mix}, 
Corollary \ref{coro:univ}, Proposition \ref{prop:tgt} and
Proposition \ref{prop:tgtrepelling} hold in this broader context.
We conclude with the following statement whose proof follows from the same arguments as on the sphere:

\begin{thm} Let $f:(\XX_1,X)\to (\XX_0,X)$ is a topologically cxc map with repellor $X$.
Assume that $X$ is locally connected and that there is an Ahlfors regular distance 
on $X$ of minimal dimension $Q$.
Then
\be
\item There exists a family of curves in $X$ of positive $Q$-modulus.
\item Thick curves are preserved by $f$.
\item There exists $\de>0$ such that every point of $X$ belongs 
to a thick curve of diameter at least $\de$.
\eb
\end{thm}


\def\cprime{$'$} \def\cprime{$'$} \def\cprime{$'$} \def\cprime{$'$}

\noindent \textsc{Peter Ha\"{\i}ssinsky,Universit\'e de Provence, CMI/LATP, 39, rue F. Joliot-Curie, 13453 Marseille cedex 13, France}
 \url{phaissin@cmi.univ-mrs.fr}
 \gap

\noindent \textsc{Kevin M. Pilgrim, Dept. Math., Indiana University, Rawles Hall, Bloomington, IN 47401 USA} \url{pilgrim@indiana.edu}

\end{document}